%Paper: alg-geom/9408011
%From: Robert Lazarsfeld <rkl@math.ucla.edu>
%Date: Wed, 31 Aug 94 13:08:45 PDT

%%%%%%%%%%%%%%%%%%%%%%%%%%%%%%%%%%%%%%%%%%%%%%%%%%%%%%%%%%%%%%%%%
%
%         LECTURES ON LINEAR SERIES
%
%                    by
%
%
%             Robert LAZARSFELD
%
%  (with the assistance of G. Fernandez del Busto)
%
%
%    Expanded notes from a course delivered at
%    the 1993 Regional Geometry Institute in
%    Park City, Utah.
%
%
%   AMS TeX 2.1
%
%
%%%%%%%%%%%%%%%%%%%%%%%%%%%%%%%%%%%%%%%%%%%%%%%%%%%%%%%%%%%%%%%

\documentstyle{amsppt}
\magnification 1100
\parindent .3in
\pagewidth{5.5in}
\pageheight{7.6 in}
\hcorrection{.2 in}
\vcorrection{.23 in}
\NoBlackBoxes

%%%%%%%%%%%%%%%%%%%%%%%% Defs of Abbreviations
%%%%%%%%%%%%%%%%%%%%%%%%%%%%%%%%%%%%%

\def \C{\bold{C}}
\def \P{\bold{P}}
\def \Q{\bold{Q}}
\def \R{\bold{R}}
\def \Z{\bold{Z}}
\def \E{\Cal{E}}
\def \F{\Cal{F}}
\def \G{\Cal{G}}
\def \I{\Cal{I}}

\def \L{\Cal{L}}

\def \A{\Cal{A}}

\def \O{\Cal{O}}
\def \eps{\epsilon}

\def \hra{\hookrightarrow}
\def \lra{\longrightarrow}

\def \ed{\quad\qed}

\def \rndup#1{\ulcorner \! #1 \! \urcorner}

\def \restr#1{\text{ } | \text{ } #1}
\def \xca#1{\bl\ni {\bf #1.}}
\def \demo#1{\sbl \ni{\bf{#1.}}}

\def \ni{\noindent}
\def \bl{\vskip 9pt}
\def \sbl{\vskip 5pt}
\def \tbl{\vskip 3pt}

%%%% Abbreviations to simulate PCMS documentstyle

\def \endxca{ }
\def \enddemo{ }
\def \head#1{\bl \ni{\bf{#1}}}
\def \endhead{ }
\def \ext{ \sbl}
\def \endext{\sbl}
\def \subhead#1{\sbl \ni{\bf{#1}}}
\def \endsubhead{ }

%%%%%%%%%%%%%%%%%%%%%%%%%%%%%%%%%%%%%%%%%%%%%%%%%%%%%%%%%%%%%%%%%%%%%%%%%%%

\centerline{\bf LECTURES ON LINEAR SERIES \footnote{Expanded
notes from a course delivered at the 1993 Regional Geometry
Institute in Park City, Utah.}}
\bl
\tbl
\centerline{\bf Robert Lazarsfeld\footnote{Partially
supported by NSF grant DMS 94-00815.}}
\bl
\centerline{ \it with the assistance of}
\sbl
\centerline{Guillermo Fern\'andez del Busto}
\bl
\bl

\head {Contents} \endhead
\sbl
\tbl \ni { Introduction}
\tbl \ni { \S 0. Notation and Conventions}
\tbl \ni { \S 1. Background and Overview}
\tbl \ni { \S 2. Reider's Theorem -- Statement and Applications}
\tbl \ni { \S 3. Building Vector Bundles}
\tbl \ni { \S 4. Reider's Theorem via Vector Bundles}
\tbl \ni { \S 5. Vanishing Theorems and Local Positivity}
\tbl \ni { \S 6. Adjoint Series and Bogomolov Instability via Vanishing}
\tbl \ni { \S 7. Algebro-Geometric Analogue of Demailly's Approach}
\tbl \ni {  References}

\sbl
\bl
\head{Introduction} \endhead
\sbl

Linear series have long played a central role in algebraic geometry,
and the basic results and techniques from their study form an
essential part of the field's culture. However the past decade has
witnessed two important new developments. First, vector bundles
have emerged as powerful tools for analyzing linear series on
curves and surfaces. These vector bundle techniques --- the most
important being Reider's method involving Bogomolov instability on
surfaces --- have led to  considerable simplifications and
extensions of classical results. A second major influence has come
from the flowering of higher dimensional geometry. One now has a
conjectural picture of how the most familiar facts about linear
series on curves should extend to arbitrary smooth projective
varieties, and some encouraging partial results have  been proved.
On the technical side, the cohomological machinery developed to
study higher dimensional varieties --- notably vanishing theorems
for $\Q$-divisors --- is proving to have applications even to
concrete questions on surfaces. These new tools promise to play an
important role in the future, but they have to a certain  extent
remained embedded in the research literature, and this has limited
somewhat their accessibility to the novice or non-expert.

\bl The purpose of these lectures is to provide a down-to-earth and
gentle introduction to some of these new ideas and techniques.
While I hope that these notes may have something to say to seasoned
geometers who wish to learn  about recent work on linear series, I
have particularly tried to gear the discussion to a novice
audience.  My intention was that with a little faith and effort,
the material here should be accessible to anyone having  finished
the standard  texts, e.g.
\cite{H2} or \cite{GH1}.
\footnote"*"{ The only pre-requisites  not covered in  \cite{H2}
that we use systematically are Chern classes of vector  bundles,
and vanishing theorems. However the latter might be taken on faith.
Griffiths and Harris \cite{GH1} provide more than enough background
for everything that appears here, except that the reader might want
to supplement their discussion of the classical Kodaira vanishing
theorem with an account (e.g. \cite{CKM}, Lecture 8, or
\cite{Kol3}, Chapters 9, 10, or \cite{EV2}) of the generalization
due to Kawamata and Viehweg.}   The underlying theme  is the search
for higher dimensional generalizations of the most basic theorems
about linear series on algebraic curves, but to keep things
elementary we work more or less entirely in dimension two.  The
philosophy is to illustrate in the setting of surfaces the various
methods that have been used to attack these questions in general,
and we end up repeatedly proving variants of one central result,
namely Reider's theorem.   Sticking to the case of surfaces allows
one to eliminate many technical complexities, and some of the
underlying ideas become particularly transparent.  I hope that
parts of the present notes might therefore provide a useful first
introduction  to the powerful and important cohomological tools of
higher dimensional geometry, as well as to an active area of
current research. To this end, I have included many exercises which
sketch further developments and applications of the material
discussed in the text. An overview of the questions we consider,
and an outline of the contents of the lectures, will be found in
\S1.

\bl
It was not my intention to produce a balanced survey of work on
linear series, and I have ignored a number of topics that might
have fit very nicely. I am particularly cogniscent of the fact
that the complex analytic side of the story is woefully
under-represented here. Starting with Demailly's ground-breaking
paper \cite{De1}, there has been an intriguing and fruitful
interplay between algebraic and analytic approaches to many of
these questions. I can only hope that the one-sided slant of the
present notes will motivate someone more able than I to give an
introductory account of the analytic viewpoint.

	\bl  I'm grateful to O. K\"uchle, V. Ma\c sek and G. Xu for helpful
suggestions, and to F. Schreyer for bailing me out of a problem in
a preliminary version of these lectures. G. Fern\'andez del Busto
--- who served as my course assistant at the Park City Institute ---
produced a preliminary version of these notes and also contributed
useful advice.  I've written an exposition on linear series once
before
\cite{L2}, but the viewpoint here is rather different. This
evolution is partly the result of contact with a number of
mathematicians to whom I owe thanks. To begin with, H. Esnault and
E. Viehweg initially explained to me a few years ago the philosophy
of how vanishing for
$\Q$-divisors could be used to produce sections of adjoint
bundles.  I've also greatly profitted from suggestions and
encouragement from J. Koll\'ar, M. Reid and Y.-T. Siu. But above
all I wish to acknowledge my debt to Lawrence Ein, with whom I have
spent the last few years working on the questions discussed in
these notes. He deserves a large share  of the credit for any
originality or utility the present lectures may possess.

\newpage
\head {\S 0. Notations and Conventions.} \endhead

\bl (0.1).  We work throughout over the complex numbers $\C$.
Varieties are assumed to be smooth and projective unless otherwise
stated.
\bl (0.2). We write $K_X$ for the canonical divisor (class) of a
smooth variety $X$. If $Z \subset X$ is a subvariety, $\I_Z \subset
\O_X$ denotes its ideal sheaf.
\bl (0.3). Given a line bundle $L$ on $X$, much of our focus will
be on the adjoint bundle $K_X + L$, for which it is traditional to
use additive notation. But then to indicate the sheaf of sections of
this bundle vanishing at a point $x \in X$, it seems most natural to
write $\O_X(K_X+ L) \otimes \I_x$. It  soon becomes difficult to
maintain notational consistency. We've finally surrendered  to a
certain amount of chaos concerning notation for line bundles,
divisors and the corresponding invertible sheaves. If $L$ is a line
bundle and $D$ is a divisor on $X$, we write $L(D)$, $L + D$ and
$O_X(L + D)$ more or less interchangably. We hope that  this will
not lead to undue confusion or annoyance. The notation $|\O_X(L +
D)|$ is reserved for the complete linear system of divisors of
sections of
$\O_X(L + D)$.  We say that a line bundle (or a divisor, or a linear
series) $L$ is {\it free} (or {\it globally generated} or {\it
basepoint-free}) at a point $x \in X$ if there is a section $s \in
H^0(X, L)$ with $s(x) \ne 0$. $L$ is {\it free} if it is free at
every point of $X$
\bl (0.4).  If $L$ is a line bundle on a smooth projective variety
$X$ of dimension $n$, we denote by $L^n$ the top self-intersection
$\int c_1(L)^n$ of the first Chern class of $L$. If $C \subset X$
is a curve in $X$, $L \cdot C$ denotes the intersection number of
$c_1(L)$ with $C$. Numerical equivalence of divisors or line
bundles is denoted by $\equiv$. Recall that $L$ is {\it nef}, or
{\it numerically effective} if $L \cdot C \ge 0$ for every
irreducible curve $C
\subset X$. By Kleiman's criterion (cf. \cite{H1}, Chapter 1), $L$
is nef if and only it is in the closure of the  cone of ample line
bundles. Note that the pull-back of a nef line bundle under any
morphism remains nef.  A nef line bundle is {\it big} if $L^n > 0$.
This is equivalent (for nef bundles) to saying that $h^0(X, kL)$
grows like $k^n$. (See \cite{CKM} or \cite{Mori} for   alternative
characterizations.)

\bl
\head {\S 1. Background and Overview}
\endhead

\sbl

The  theme of these lectures is the search for higher dimensional
generalizations of the most familiar and elementary facts about linear
series on curves. In this section we introduce the basic questions, and
give a brief overview of their history and current status. The detailed
contents of the notes are summarized at the end of the section.
\sbl

To set the stage, we start by recalling the story in dimension one:
\proclaim{Theorem 1.1} Let  $C$ be a compact Riemann surface of  genus $g$.
\sbl  (A). If $g\geq 2$ then the canonical bundle $K_C$ is globally
generated,  and  the pluri-canonical series $|mK_C|$ are very ample when $m
\ge 3$.
\sbl (B). Let $N$ be a line bundle on $C$, with $\deg (N)=d$.  If $d\geq
2g$ then $N$ is  globally generated, and if $d\geq 2g+1$ then $N$ is very
ample.
\endproclaim
\bl
\demo{Proof} We focus on statement (B). Suppose that $d\geq 2g$,  and let
$P\in C$ be a fixed point. We need to show that there exists a section
$s\in \Gamma (C,N)$  with $s(P)\neq 0$. Consider to this end the exact
sequence of sheaves
$$
0\lra N(-P) \lra N  \lra N \otimes \O_P  \lra 0 .
$$ It induces
in cohomology the exact sequence
$$H^0(C,N) \overset {e_P} \to \lra
  H^0(C,N \otimes \O_P ) \lra H^1(C,N(-P)) \phantom{.},$$
with $e_P$ evaluation at $P$. Now $H^0(C,N \otimes \O_P) \cong \C$, and it
suffices  to show that $e_P$ is surjective. But this follows from the fact
that $\deg N(-P) \geq 2g-1$, since in this case we have that
$H^1(C,N(-P))=0$. The proof of very ampleness is similar, as is statement
(A). (cf. \cite{H2}, IV.3.2). \qed \enddemo

\bl

The attempt to generalize this theorem to higher dimensions is  a very
fundamental  and interesting problem. In setting (A), it is relatively
clear what to look for. The higher dimensional analogue of a curve of genus
$\ge 2$ is a variety of general type. Taking into account the expectation
that one might want to work  on a particularly tractable birational model,
we can state the question as
\proclaim{Problem A } Study pluri-canonical series $|mK_X|$ on ``nice"
varieties $X$  of general  type. \endproclaim

\bl By contrast, the higher dimensional analogue of statement (B) has only
recently come into focus. The key is to rephrase (1.1)(B) without
explicitly bringing in the genus of the curve.  To this end, note that if
$N$ is a line bundle of degree
$2g$ on a curve $C$ of genus $g$, then we can write $N$ in the form
$$ N = K_C + 2A$$ for some ample line bundle $A$ on $C$. Similarly, if
deg$(N) = 2g+1$, then $N = K_C + 3A$. Hence  the natural question
generalizing statement (B) is
\proclaim{Problem B } Let $X$ be a smooth projective variety of dimension
$n$. Study the {\it adjoint linear series} $$| K_X+L|$$  where $L$ is a
suitably positive line bundle  on $X$.
\endproclaim

\bl \ni For example, given an ample line bundle $A$ on $X$, one might take
$L=(n+1)A$ or $L=(n+2)A$. As a start, one would like to understand when the
adjoint bundles in question are  globally generated or very ample.

\bl Pluricanonical mappings of surfaces of general type were studied by
Kodaira \cite{Kod} and Bombieri  \cite{Bomb} in the late 1960's and early
1970's.  In 1988, many of their results emerged as special cases of a
theorem of Reider \cite{Rdr} concerning adjoint linear series. Precise
statements  appear in \S2. Suffice it to say here that  Problems A  and B
are by now quite well understood on surfaces.

\bl In higher dimensions, naturally enough, much less is known. But in
recent years two conjectures of Fujita have attracted a great deal of
interest:

\proclaim{Fujita's Conjectures } \cite{Fuj1} Let $X$ be a smooth projective
variety of dimension $n$.
\vskip 3pt (A). Assume that $X$ is minimal and of general type, i.e.
suppose that   the canonical bundle $K_X$ is nef and big. Then
$\O_X(mK_X)$ is globally generated for $m \ge n+2$.
\vskip 4pt (B). Let $A$ be an ample line bundle on $X$ (which is now not
assumed to be minimal or of general type). Then $K_X + (n+1)A$ is globally
generated, and $K_X + (n+2)A$ is very ample.
\endproclaim
\bl \ni For curves this is essentially the content of Theorem 1.1, and on
surfaces these statements follow from Reider's theorem. But in general
Fujita's conjectures remain open as of this writing. Existing results are of
two sorts. First, there are effective statements in all dimensions due to
Demailly \cite{De1}, Koll\'ar \cite{Kol2} and most recently Siu
\cite{Siu2},  which however are exponential in the dimension $n$. In another
direction, one can stick to the simplest case of global generation on
threefolds, and ask for statements closer to the bounds predicted by Fujita.
Theorems of this sort are given by Ein and Lazarsfeld in
\cite{EL2}, with Ma\c sek in \cite{ELM}, and by Fujita in
\cite{Fuj2}. We refer to Ein's paper \cite{E} for an overview of some of this
work. Siu's paper on Matsusaka's theorem \cite{Siu1}
develops some related ideas  (cf. Exercise 7.7).

\sbl It is probably well to stress from the outset why it is that Problems
A and B, while completely elementary on curves, become  more subtle in
higher dimensions. The proof of  Theorem 1.1 reduces the
question to the vanishing of a cohomology group, and this goes through in
the general setting without any problem:
\proclaim {Proposition 1.2}   Let $L$ be an ample line bundle on a smooth
projective variety $X$ (of arbitrary dimension). Then $K_X+L$ is globally
generated if and only if
$$H^1(X,\O_X(K_X+L)\otimes  \I_x)=0$$ for all points $x\in X$, where
$\I_x$ is the ideal sheaf of   $\{ x\}$.
\qed
 \endproclaim
\sbl \ni (There is of course an analogous criterion for $K_X + L$ to be
very ample.) If $X$ is  a curve, then $\O_X(K_X+L)\otimes   \I_x$ is
locally free, and one deduces the required vanishing from general facts
about the cohomology of line bundles. However when $X$ is a variety of
dimension $\geq 2$, the sheaf $\O_X(K_X+L)\otimes \I_x$  is no longer
invertible and it becomes considerably harder to control the group
$H^1(X,\O_X(K_X+L)\otimes \I_x)$. So in higher dimensions the issue in effect
is to prove a Kodaira-type vanishing theorem for certain non-invertible
sheaves.

\bl The ``classical" approach to these questions on surfaces   is to apply
vanishing theorems involving numerical connectedness on a blow-up of
$X$ to control the required cohomology groups. This  is discussed in Reid's
lectures. During the 1980's    some new  techniques emerged for studying
linear series on curves and surfaces, revolving around the geometry of
vector bundles. Methods along these lines were used for example to study
projective normality of curves in \cite{GL}, and the geometry of curves on
K3 surfaces in \cite{L1}. But certainly the most important
development  was Reider's application \cite{Rdr} of Bogomolov's instability
theorem to study adjoint series on surfaces. Reider's method -- which has
its antecedents in a proof by Mumford of Kodaira vanishing
\cite{Reid} -- has largely superceeded the classical approach using
numerical connectedness \`a la Ramanujam. We remark that a number of
geometers have  attempted to  apply vector bundles in a similar manner to
study linear series on varieties of dimension three and higher, but so far
these efforts have   not met with  success.

\bl At about the same time that vector bundle techniques were  developed to
study linear series on surfaces,  Kawamata, Reid and Shokurov (among
others) introduced some powerful but subtle cohomological techniques to
analyze the asymptotic behavior of pluricanonical and other linear series
on varieties of all dimensions. These techniques --  which we shall refer
to as the KRS package -- are based on the Kawamata-Viewheg vanishing
theorem for $\Q$-divisors, and form one of the central tools of
contemporary higher dimensional geometry.  A typical result here is a
theorem of Kawamata and Shokurov to the effect that if $X$ is a smooth
minimal variety of general type, then $|mK_X|$ is free for all $m \gg 0$.
We refer to \cite{CKM} or \cite{KMM} for references and an overview of this
circle of ideas.  Ein and the  author remarked that one could use the KRS
machine to recover (some of) Reider's results. This in turn opened the door
to the higher dimensional theorems appearing for instance in \cite{Kol2}
and \cite{EL2}.

\bl The goal of these lectures is to explain some of the new methods that
have been introduced to study Problems A and B on surfaces and higher
dimensional varieties. To keep things  elementary, we will work more or
less exclusively on surfaces, and the plan is to discuss one central result
(viz. Reider's theorem) from many points of view.  To begin with,  the
vector bundle techniques apply here, and we go through these in some
detail. However we also present in the surface setting some of the  methods
developed for  higher dimensions. Thus we devote considerable attention to the
circle of ideas surrounding vanishing theorems for $\Q$-divisors.  We also
explain -- still in the case of surfaces -- some algebro-geometric analogues of
Demailly's analytic approach, as well as some local results from \cite{EL3} and
\cite{EKL}.

\bl

The detailed contents of these notes is as follows. In \S2 we give the
statement of Reider's theorem, and present some of its applications. The
following two sections are devoted to the proof of this result via vector
bundles: \S3 discusses Serre's method (and related techniques) for
constructing vector bundles, and in \S4 Bogomolov's instability theorem is
applied to prove Reider's statement. The focus then turns to vanishing
theorems. We consider in  \S5 the questions involving Seshadri constants
and local positivity that arise if one tries to apply vanishing theorems in
the most naive possible way to study adjoint series. Section 6 revolves around
vanishing theorems for $\Q$-divisors and the Kawamata-Reid-Shokurov
machine: we use these techniques to reprove (parts of) Reider's theorem, as
well as the rank two case of Bogomolov's theorem. Finally,  \S7 is devoted
to the algebro-geometric analogues of Demailly's approach to these
questions. Further applications and extensions of the material are outlined
in numerous exercises scattered throughout the notes.

\bl \bl
\head{\S2. Reider's Theorem -- Statement and Applications}
\endhead

\bl

In this section, we state Reider's theorem \cite{Rdr}, and present some
simple applications. Reider's result forms the core of the present notes:
most of the subsequent sections are devoted to various different proofs
of the theorem or its corollaries. Here we try to convey some feeling for
its power and scope.

 \bl
The theorem in question gives a very precise geometric explanation for the
failure of an adjoint bundle on a surface to be globally generated or
very ample.

\proclaim{Theorem 2.1 \rm (\cite{Rdr})} Let $X$ be a smooth projective
surface, and let $L$ be a nef line bundle on $X$.
\sbl
\ni (i). Assume that $ L^2 \ge 5$, and that the adjoint series $|K_X +
L|$ has a base-point at $x \in X$. Then there exists an effective divisor
$D \subset X$ passing through $x$ such that either
 $$\aligned  D \cdot L = 0 \quad &\text{and} \quad D^2 = -1; \ \ \text{
or} \\
	D \cdot L = 1 \quad &\text{and} \quad D^2 = 0. \endaligned \tag 2.2
$$
\ni (ii). If $L^2 \ge 10$, and if $x, \ y \in X$ are two points (possibly
infinitely near) which fail to be separated by $|K_X + L|$, then there
exists an effective divisor $D \subset X$ through $x$ and $y$ such that
either
$$\aligned D \cdot L = 0 \quad &\text{and} \quad D^2 = -1 \ \text{ or } \
-2; \
\ \text{or} \\
	D \cdot L = 1 \quad &\text{and} \quad D^2 = 0 \ \text{ or } \ -1; \ \
\text{or}  \\
	D \cdot L = 2 \quad &\text{and} \quad D^2 = 0 \endaligned \tag 2.3
$$
\endproclaim
\ni Note that we do not assume that the base point $x$ in (i) is
isolated. Hence subject to the numerical inequalities on the nef line
bundle $L$, the statement applies as soon as $\O_X(K_X + L)$ fails to be
globally generated. Similarly, (ii) serves as a criterion for the bundle
in question to be very ample.

\demo{Example 2.4} Typical examples of a divisor satisfying (2.2) may be
obtained by taking $D$ to be an exceptional curve of the first kind (in
the first case), or as the fibre of a ruled surface (in the second).
\enddemo

\bl
\xca{Exercise  2.5}  (i). Show that if $D$ is one of the divisors
satisfying (2.2) described in Example 2.4, then $D$ is necessarily in the
base locus of
$|K_X + L|$.
\sbl
 (ii).  Find  examples of  divisors $D$ satisfying the various
possibilities in (2.3), for which $\O_X(K_X + L)$ indeed fails to be very
ample.
\sbl
 (iii). Show that the numerical hypothesis on $L^2$ in  Theorem 2.1
cannot  in general be weakened. \qed

\endxca

\bl Reider's theorem leads to a simple numerical criterion for an adjoint
bundle to be globally generated or very ample:

\proclaim{Corollary 2.6} Let $L$ be an ample line bundle on  a smooth
projective surface $X$. Assume that $L^2 \ge 5$, and that $L \cdot C \ge
2$ for all irreducible curves $C \subset X$. Then $\O_X(K_X+L)$ is
globally generated. If
$L^2\geq 10$ and  $L \cdot C\geq 3$ for all  $C \subset X$, then
$\O_X(K_X+L)$ is very ample. \ed \endproclaim

\bl
\ni This in turn implies the dimension $n = 2$ case of Fujita's
conjecture:

\proclaim{Corollary 2.7} If $A$ is an ample line bundle on $X$, then
$|K_X+3A|$ is free, and $|K_X+4A|$ is very ample.
\endproclaim
\demo {Proof} In fact, take $L = 3A$. Then $L^2 = 9 A^2 \ge 9$, and $L
\cdot C = 3 A \cdot C \ge 3$ for any curve $C \subset X$. So Corollary
2.6 shows that
$\O_X(K_X + L)$ is globally generated, and similarly $\O_X(K_X +4A)$ is
very ample.  \qed \enddemo

\bl As we indicated in the previous section,   attempts to generalize
Reider's theorem to higher dimensions have mainly focused on extending
the statements of Corollaries 2.6 and 2.7. At the moment, one doesn't
even have any clear conjectures as to what might be the analogue of the
more precise  information contained in Theorem 2.1. However in the case
of surfaces, these restrictions on the self-intersection of  $D$ are very
powerful. Indeed,   we will now see that they allow one to deduce from
Reider's theorem many of the classical facts concerning pluricanonical
models of surfaces of general type.

\bl Let $X$ be a  surface of general type. It is very natural and
important to try to understand the geometry of the pluricanonical
rational mappings
$$
\Phi_m = \Phi_{|mK|} : X  \dashrightarrow \P = \P H^0(mK_X)
$$ defined by multiples of the canonical bundle $K_X$. Here Reider's
theorem leads to an extremely quick proof of some of the fundamental
results of Kodaira and Bombieri:

\proclaim{ Corolllary 2.8} {\rm (\cite{Kod}, \cite{Bomb})} Assume that
$X$ is minimal i.e. not the blowing up of some other smooth surface at a
point. Then:
\bl
\ni (i). The  bundle $\O_X(mK_X)$ is globally generated (i.e. $\Phi_m$ is
a morphism) if $m \ge 4$, or if $m \ge 3$ and $K_X^2 \ge 2$.
\bl
\ni (ii). $\Phi_m$ is an embedding away from $(-2)$-curves if $m \ge 5$,
or if $m
\ge 4$ and $K_X^2 \ge 2$, or if $m \ge 3$ and $K_X^2 \ge 3$. \endproclaim

\ni \demo{Remark}  By an embedding away from $(-2)$-curves, we mean   a
morphism which is one-to-one and unramified away from a divisor
$Z\subseteq X$ consisting of smooth rational curves with
self-intersection $-2$.    Catanese and Reider
\cite{C} have used Reider's method to prove the more precise theorem of
Bombieri that $K_V$ is very ample on  the canonical model
$V$ of $X$  obtained by blowing down all the $A-D-E$ cycles $Z$ to
rational double points. \enddemo

\demo{Proof of Corollary 2.8} We will only consider the global
generation of
$|mK_X|$, leaving the proof of the second assertion as an exercise.
\sbl Recall that on a minimal surface of general type $X$, the canonical
bundle
$K_X$ is nef \cite{BPV}, III.2.3.  Hence Theorem 2.1 applies  to
$L=(m-1)K_X$. The numerical hypotheses of (2.1)(i)  are satisfied thanks
to the conditions on
$m$ and $K_X^2$. Thus if $|mK_X|$ has a base point, then there exists an
effective divisor $D\subseteq X$ such that either $(m-1)K_X \cdot D = 1$
or else
$$(m-1)K_X \cdot D = 0 \  \text{and} \ D^2 = -1.$$ The first possibilty is
excluded by the assumption that $(m-1) \ge 2$. As for the second, if
$K_X\cdot D=0$ then
 $D^2=D\cdot (D+K_X)\equiv 0 \ \text{ (mod 2) }$ by adjunction,
contradicting
$D^2=-1$. \ed \enddemo
\bl We refer  to \cite{Rdr} for further applications of Reider's method to
pluricanonical mappings. The following exercises present some other
applications of Theorem 2.1, and the conference proceedings \cite{SBL}
contain a sampling of some more recent developments.

\bl
\xca{Exercise 2.9}  ({\bf Embeddings of Abelian Surfaces.}) Prove the
following result of Ramanan \cite{Ram}. Let $X$ be an abelian surface
which contains no elliptic curves. If $L$ is an ample  line bundle on $X$
such that $L^2 \ge 10$, then $L$ is very ample.  (It follows for example
that on a sufficiently general abelian surface, a polarization of type
$(1, d)$ is very ample if $d \ge 5$. In particular, taking d = 5, there
exist smooth abelian surfaces in $\P^4$.) See
\cite{LB} Chapt. 10, \S4, for more precise statements, and the following
exercise for a generalization. The corresponding questions on higher
dimensional abelian varieties are considered in \cite{BLR} and
\cite{DHS}. \qed
\endxca

\bl
\xca{Exercise 2.10} ({\bf Linear Series on Minimal Surfaces of Kodaira
Dimension  Zero.})  In this exercise,  $X$ denotes a smooth projective
surface whose canonical bundle  $K_X$ is numerically trivial. This
hypothesis is satisfied e.g. by Abelian, K3 and Enriques surfaces; see
\cite{BPV}, IV.1, for the complete list.  Fix an ample line bundle $L$ on
$X$.
\sbl
 (i).  Assume that $L^2 \ge 5$. Prove that $\O_X(K_X + L)$ fails to be
globally generated if and only if there exists an irreducible curve $E
\subset X$ with $p_a(E) = 1$ and $E \cdot L = 1$.
\sbl
 (ii).  Assume that $L^2 \ge 10$, and that $|K_X + L|$ is free. Prove that
$\O_X(K_X + L)$ fails to be very ample if and only if there exists a
reduced curve
$E \subset X$, with $p_a(E) = 1$ and $E \cdot L = 2$. [Note that
the intersection form  on $X$ is even, i.e. $D^2$ is even for every
divisor $D$ on $X$.]
\sbl
\ni See \cite{Rdr}, Proposition 5, for a more precise statement due to
Beauville.
 \qed \endxca

\bl

We conclude this section with some interesting open problems of an
algebraic nature that arise in connection with Reider's theorem.
Returning for an instant to the one-dimensional case, let $N$ be a line
bundle of degree $d$ on a smooth projective curve $C$ of genus $g$. A
classical theorem of Castelnuovo, Mattuck and Mumford  asserts that if $d
\ge 2g+1$, then $N$ is normally generated, i.e. the natural maps
$Sym^m(H^0(N)) \lra H^0(mN)$ are surjective for all $m \ge 0$.
Furthermore, if $d \ge 2g+2$, then in the embedding $C \subset
\P{\text{\rm H}}^0(N)$ defined by $N$, the homogeneous ideal
$I_C$ of $C$ is generated by quadrics. The famous theorems of Noether and
Petri (cf. \cite{Mfd}, \cite{ACGH} Chapter III, \cite{GL}) give
analogous statements for the canonical bundle, and Green \cite{Grn} has
shown that at least conjecturally the whole picture extends to higher
syzygies as well. (See also Exercise 3.5.)

\bl The question then arises whether similar statements hold for adjoint
and pluricanonical bundles on surfaces and higher dimensional varieties.
In the two dimensional case, the natural thing to hope for here is the
following:
\proclaim{Conjecture 2.11} {\rom{({\bf Mukai}.)}} Let $A$ be an ample line
bundle on a smooth projective surface $X$, and let $P$ be any nef line
bundle. Then
$\O_X(K_X + 4A + P)$ is normally generated, and in the embedding $X
\subset \P$ defined by $\O_X(K_X + 5A + P)$, the homogeneous ideal $I_X$
of $X$ is generated by quadrics.
\endproclaim
\bl
\ni One would also like results in the spirit of (2.1) dealing with the
adjoint bundles $K_X + L$ subject only to numerical conditions on $L$.
For instance, one might hope that if $L$ satisfies the hypotheses of
(2.1)(ii), then
$\O_X(K_X + L)$ is normally generated. There should also be analogous
statements for higher syzygies.

\bl When $X$ is a ruled surface, Butler \cite{But} proves that $K_X + 5A +
P$ is projectively normal, and he also obtains results  for generation by
quadrics and higher syzygies. For ``hyper-adjoint" bundles of the form $K
+ mB + P$ where
$B$ is {\it very ample}, results in arbitrary dimension are given  by Ein
and Lazarsfeld in \cite{EL1}. In fact, if $V$ is a smooth projective
variety of dimension $n$, then $K_V + (n+1)B + P$ is normally generated,
and $K_V + (n+2)B + P$ defines an embedding in which the homogeneous
ideal of $V$ is defined by quadrics. (Compare Exercise 5.15. \cite{EL1}
also contains analogous statements for higher syzygies.)  However it
seems that new ideas will be needed to tackle Conjecture 2.11.  It would
be already very interesting to get a result along the lines of (2.11) but
with weaker numbers.

 \sbl
\bl
\head{\S 3. Building Vector Bundles}\endhead
\sbl

This section is devoted to preparations for the first proof of Reider's
theorem, via vector bundles. The strategy for using vector bundles to study
linear series involves two steps:
\sbl

\item {$\bullet$} Encode the geometric data at hand in a vector bundle $E$; and
\vskip 3pt
\item{ $\bullet$} Study the geometry of $E$
 (e.g. stability, or sub-bundles, or endomorphisms) to arrive at  the desired
conclusions.  \sbl

\ni Here we focus on the first point, and discuss techniques for constructing
bundles. Specifically, after some warm-up with extentions of line bundles,  we
study a method introduced by Serre, which underlies Reider's argument. In the
exercises, we sketch several other constructions and applications. A good
reference for the general homological machinery appearing here is \cite{GH1},
Chapter 5, \S3 and \S4.

\sbl
\bl
%\ni {\it Extensions of Line Bundles}
\subhead{Extensions of Line Bundles.}\endsubhead
\sbl

To set the stage, we begin with the simplest technique for producing vector
bundles, namely as extensions of invertible sheaves.  Let $X$ be an irreducible
projective variety, and let $L$ and $M$ be line bundles on
$X$. Recall that an {\it extension} of $L$ by $M$ is a short  exact sequence of
sheaves:
$$
\CD 0 @>>> M @>>> E @>>> L @>>> 0,
\endCD \tag 3.1$$ so that $E$ is a rank two vector bundle on $X$. Two such
extentions are {\it  equivalent} if there is a map between them inducing the
identity on the outer terms. The set of equivalence classes of such extensions
can be given the structure of a complex vector space, denoted $\text{\rm
Ext}^1(L,M)$. The zero element of $\text{\rm Ext}^1(L,M)$ corresponds to the
split sequence. (Analogous considerations apply when $L$ and $M$
are arbitrary coherent sheaves, although of course in this case
the middle term in (3.1) isn't in general locally free.)

\bl

Assuming as we are that $L$ and $M$ are line bundles, one has the basic
isomorphism:
$$
\text{\rm Ext}^1(L,M) = H^1(X, L^* \otimes M),  \tag 3.2
$$ whose proof we outline in Exercise 3.3. For the present purposes, the
importance of (3.2) is that the bundle $E$ serves as a geometric realization of
a class $e \in H^1(X, L^* \otimes M)$. This is illustrated in Exercise 3.5,
where we indicate how ideas along these lines may be used to prove Noether's
theorem and various generalizations concerning the projective normality of
algebraic curves.

\bl
\xca{Exercise 3.3} ({\bf Extension Classes } -- Compare \cite{H2} Ex. III.6.1,
and \cite{GH1}, pp. 722 - 725.) We keep the notation just introduced. In
particular, $L$ and $M$ are locally free sheaves of rank one on the projective
variety $X$.

\sbl
 (i).  We start by defining the map
$$
 \text{\rm Ext}^1(L, M) \lra H^1(X, L^* \otimes M) \tag 3.3.1$$ required for
(3.2). To this end, fix an extension  of $L$ by $M$. Prove that there exists a
covering  of
$X$  by open sets $\{ U_i \}$ on which the restriction of (3.1) splits.  Show
that on the intersections $U_i \cap U_j$, the difference of two such local
splittings  determines  homomorphisms
$L(U_i \cap U_j)\lra M(U_i \cap U_j)$, which in turn yields a \v Cech cocycle
in $ Z^1( \{U_i
\}, L^*
\otimes M)$. Prove that the resulting cohomology class  $e
\in H^1(X, L^* \otimes M)$ --- which is called the {\it extension class} of the
extension (3.1) --- is independent of the choices made, and vanishes if and
only if (3.1) splits. This defines (3.3.1).

\sbl
 (ii).  Show that the map (3.3.1) is an isomorphism. [It may be helpful to
remark that one can take transition matrices for the bundle $E$ appearing in
(3.1) to be triangular, with the cocycle representing the extension class  as
the off-diagonal entry.]

\sbl
 (iii). Prove that map (3.3.1) can alternatively be defined as follows. Given
(3.1), tensor through by $L^*$ to get $0 \lra L^* \otimes M  \lra  L^*
\otimes E \lra \O_X \lra 0$. Then the corresponding extension class is the
image of the constant function $1_X \in H^0(X, \O_X)$ under the connecting
homomorphism $\delta : H^0(X, \O_X) \lra H^1(X , L^* \otimes M)$.

\sbl
 (iv). Note that cup product determines a map $H^0(X, L) \otimes H^1(X, L^*
\otimes M) \lra H^1(X,M)$, or equivalently a homomorphism
$$H^1(X, L^* \otimes M) \lra Hom( H^0(X, L) ,  H^1(X, M)) .\tag 3.3.2$$ Verify
that  (3.3.2) may be interpreted as the map which takes an extension to the
connecting homomorphism it determines.  \qed \endxca

\bl
\xca{Exercise 3.4}  Let $E$ be a rank two vector bundle on $\P^n$ which admits
a
nowhere vanishing section. Prove that if $n \ge 2$, then $E$ is a direct sum of
line bundles. (The same statement is true when $n = 1$, but the proof is a
little more delicate. In fact, a theorem of Grothendieck (cf. \cite{OSS},
I.2.1) states that any vector bundle on $\P^1$ splits as a direct sum of line
bundles.) \qed
\endxca

\bl
\xca{Exercise 3.5} ({\bf Noether's Theorem and Generalizations}.) Let $C$ be a
smooth projective curve of genus $g \ge 2$. A classical theorem of Noether
states that if
$C$ is non-hyperelliptic, then the natural maps
$$\rho_m : Sym^m(H^0(C,K_C)) \lra H^0(C, mK_C)$$ are surjective for $m \ge 2$.
In other words, every pluri-canonical differential form can be expressed as a
polynomial in  one-forms, or equivalently,  $C$ is projectively normal in its
canonical embeding $C\subset \P^{g-1}$. In this exercise we will outline a
proof
of this statement via vector bundles. There are quicker approaches to Noether's
theorem, but with only a little extra effort the present argument yields also
some substantial generalizations of the classical results (see (v) below).
These appear in \cite{GL}, from which this exercise is adapted. For simplicity
we will focus on the surjectivity of $\rho_2$; the case of
$\rho_m$ for $m > 2$ is similar but more elementary.

\sbl
 (i). Let $C$ be any smooth curve of genus $g$, and assume that $\rho_2$ is not
surjective. Then the map
$ H^0(C, \O(2K_C))^* \lra H^0(C, \O(K_C))^* \otimes H^0(C, \O(K_C))^*$ dual to
multiplication has a non-trivial kernel. Using Serre duality and Exercise
3.3.(iv), show that there exists a non-split extension
$$0 \lra \O_C \lra  E \lra \O_C(K_C) \lra 0 \tag 3.5.1$$ which is exact on
global sections. In particular, $h^0(C, E) = g+1$. The plan is to show that the
existence of such a bundle forces $C$ to be hyperelliptic.
\sbl
 (ii). Now fix a point $P \in C$. Show that there exists a non-zero section $s
\in H^0(C, E)$ vanishing at $P$. Denoting by $D$ the effective divisor
(containing $P$) on which
$s$ vanishes, one has an exact sequence:
$$ 0 \lra \O_C(D) \lra E \lra \O_C(K_C - D) \lra 0. \tag 3.5.2
$$ Using the fact that (3.5.1) is non-split, show that $\O_C(D) \not \cong
\O_C(K_C)$.
\sbl
 (iii). Compare the estimates on $h^0(C, E)$ obtained from (3.5.1) and (3.5.2)
to prove that $deg(D) \le 2 r(D)$, where as usual $r(D) = h^0(C, \O_C(D)) - 1$.
Finally, conclude from Clifford's theorem that $C$ is hyperelliptic. This
completes the proof of Noether's theorem.
\sbl (iv). Prove the theorem of Castelnouvo, Mattuck and Mumford that if  $L$
is a line bundle of degree $\ge 2g + 1$ on a curve $C$ of genus $g$, then $L$
is normally generated, i.e. the maps
$Sym^m(H^0(C,L)) \lra H^0(C, mL)$ are surjective for $m \ge 1$.
\sbl
 (v). A result going back to Segre (cf. \cite{MS}) states that if $F$ is a
rank two vector bundle of degree $d$ on $C$, then $F$ contains a line
sub-bundle $A
\subset F$ with
$ deg(A) \ge [\frac{d - g + 1}{2}]$. Using this, prove the theorem of \cite{GL}
that if $L$ is a very ample line bundle with
$$ deg(L) \ge 2g + 1 - 2h^1(L) - \text{Cliff}(C), \tag 3.5.3
$$ then $L$ is normally generated. Here
$$\text{ Cliff}(C) = \min \left \{ deg(A) - 2r(A) \ | \ h^0(A) \ge 2, deg(A)
\le g-1
\right \}$$ denotes the Clifford index of $C$, which  measures how general $C$
is from the point of view of moduli. (For instance, Cliff$(C) = 0$ if and only
if $C$ is hyperelliptic; Cliff$(C) = 1$ iff $C$ is either trigonal or a smooth
plane quintic; and if $C$ is a general curve of genus $g$, then Cliff$(C) = [
\frac{g-1}{2}]$.) Note that (3.5.3) contains the theorems of Noether and
Castelnuovo-Mattuck-Mumford as special cases. \qed

\endxca

\sbl
\bl

%\ni{\it The Serre Construction}
\subhead{The Serre Construction.} \endsubhead

\sbl We henceforth assume that $X$ is a smooth projective surface. The bundles
on $X$ that arise as extensions of line bundles are rather special, as Exercise
3.4 suggests.  Reider's theorem requires a more general construction. It was
introduced by Serre, and  applied notably in the analysis of codimension two
subvarieties of projective space (cf. \cite{OSS}, I.5). The case of surfaces,
which has a somewhat different flavor, was studied by Griffiths and Harris
\cite{GH2}.

\bl Consider to begin with a rank two vector bundle $E$ on $X$, with
$\text{det}(E)= L$. Suppose that
$s \in \Gamma(X, E)$ is a section of $E$ that vanishes on a finite set. Denote
by
$$Z = Z(s) \subset X$$ the zero-scheme of $s$:  locally one may view $s$  as
given by a vector $s = (s_1, s_2)$ of regular functions, and $Z$ is defined
locally by the vanishing of $s_1$  and $s_2$. Remark that  the Koszul complex
associated to $(E, s)$ determines a short exact sequence
$$0 \lra \O_X \overset{\cdot s}\to \lra E \lra L \otimes \I_Z \lra 0 \tag 3.6$$
resolving $L \otimes \I_Z$, where $\I_Z$ is the ideal sheaf of $Z$. Somewhat
abusively, we will refer to (3.6) as the Koszul complex arising from $s$.

\bl Serre addresses the possibility of making the inverse construction.
Specifically, fix a line bundle $L$ on $X$, plus a finite scheme $Z \subset X$.
We pose
\proclaim{Question 3.7}  When does there exist a rank two vector bundle $E$ on
$X$, with
$\text{det}(E) = L$, plus a section $s \in \Gamma(X, E)$ such that $Z(s) = Z$?
\endproclaim
\bl
\ni Clearly there are local obstructions to finding $E$: for instance, $Z$ must
be locally defined by two equations. But there are also some very interesting
global conditions, and these are the key to Reider's method.

\bl Serre's idea is that one should try to construct $E$ via (3.6), as an
extension of
$L \otimes \I_Z$ by $\O_X$.  As above, we  consider the group
$$\text{\rm Ext}^1(L \otimes \I_Z , \O_X)$$ parametrizing all such. Given $e
\in
\text{\rm Ext}^1(L \otimes \I_Z , \O_X)$ we denote by
$\F_e$ the sheaf arising from the corresponding extension:
$$ 0 \lra \O_X \ \lra \F_e \lra L \otimes \I_Z \lra 0. \tag 3.8
$$ Thus $\F_e$ is a torsion-free $\O_X$-module of rank two. Suppose it happens
that
$\F_e$ is actually locally free for some $e \in \text{\rm Ext}^1(L \otimes
\I_Z ,\O_X)$. Then the  map $\O_X \lra \F_e$ determines a section $s \in
\Gamma(X,\F_e)$ with (3.8) as the corresponding Koszul complex. In particular,
$Z(s) = Z$. Therefore $(\F_e, s)$ gives the required solution to (3.7).  Thus
the essential point is to determine when there exists an element $e \in
\text{\rm Ext}^1(L
\otimes \I_Z , \O_X)$ such that $\F_e$ is a vector bundle.

\bl We begin with a criterion for the failure of $\F_e$ to be locally free.
\proclaim{Proposition 3.9} Given an element $e \in {\text{\rm Ext}^1} (L
\otimes \I _Z, \O _X)$, the  correponding sheaf $\F_e$ fails to be locally free
if and only if there exists a proper subscheme
$Z^{\prime}  \subsetneqq  Z$ (possibly $Z^{\prime}=\emptyset$)
such that
$$e\in \text{Im} \left \{ \text{{\rm Ext}}^1( L \otimes
\I_{Z^{\prime}}, \O _X)
\lra  \text{{\rm Ext}}^1(L \otimes \I _Z, \O _X) \right \} . \tag
3.9.1
$$
\endproclaim

\bl
\demo{Remark 3.10} Let us explicate the map in (3.9.1). If $Z^{\prime}\subseteq
Z$, then
$L \otimes \I _Z \subseteq  L \otimes \I_{Z^{\prime}}.$  Thus
starting with an extension of $L \otimes \I _{Z'}$ by $\O _X$,
one can pull it back in an evident way to get an extension of $L
\otimes \I _Z$ by $\O _X$. This gives rise to the homomorphism
appearing in (3.9.1). \enddemo
\bl

We now give a proof of Proposition 3.9 drawing on some facts from local
algebra. A more geometric approach to the main implication is sketched in
Exercise 3.19.
\demo{Proof of Proposition 3.9}  We start with some general
remarks of a homological nature, referring to \cite{OSS}, II.1.1 (especially
pp.
148-149), for details and proofs. Let $\F$ be a torsion-free sheaf on the
smooth surface
$X$. Then the set of points at which $\F$ fails to be locally free is finite or
empty. Denote by
$\F^{**}$ the double dual of $\F$. Then $\F^{**}$ is reflexive, hence a vector
bundle. (In general, a reflexive sheaf on a smooth variety is locally free off
a set of codimension $\ge 3$.) Furthermore, the natural injection
$\mu : \F \lra \F^{**}$ fails to be an isomorphism exactly over the set where
$\F$ is not locally free. It is suggestive to think of $\F^{**}$ as a sort of
``canonical desingularization" of $\F$.
\bl
 Turning now to the situation of the Proposition, write $\F =\F_e$, and assume
that $\F$ fails to be locally free.  The extension defining $\F$, plus the
map $\mu : \F \lra \F^{**}$, give rise to an exact commutative diagram:
$$
\CD @.     @.     0    @.    0 \\ @.     @.     @VVV   @VVV    \\ 0 @>>> \O_X
@>>>
\F @>>> L \otimes \I_Z @>>> 0 \\ @.     @|     @V{\mu}VV  @VVV    \\ 0 @>>>
\O_X  @>>> \F^{**} @>>> L \otimes \I_{Z^{\prime}} @>>> 0  \ . \\ @. @. @VVV
@VVV \\ @. @.
\tau @= \tau \\ @.     @.     @VVV   @VVV    \\ @.     @.     0    @.    0
\endCD \tag 3.11
$$ Here $\tau = \F^{**} / \F$ is a finite sheaf supported on the set where $\F$
fails to be locally free. The map $\O_X \lra \F^{**}$ determines a section
$s^\prime
\in \Gamma(X ,\F^{**})$ vanishing on a finite scheme $Z^{\prime} =_{\text{def}}
Z(s^\prime)$, and the second row of (3.11) is  the corresponding Koszul
complex. It follows from Remark 3.10  that if
$e^{\prime} \in \text{Ext}^1(L \otimes I_{Z^\prime}, \O_X)$ denotes the
extension  class of the middle row, then $e^{\prime}$ maps to the given
extension class $e \in
\text{Ext}^1(L \otimes I_{Z}, \O_X)$. So to prove the first
implication in Proposition, it suffices to show that  $Z^{\prime}$ is a proper
subscheme of $Z$. But $\tau = \text{ coker} (\mu) \ne 0$ since $\F$ is not
locally free, and the  right hand column of (3.11) then shows that $Z^{\prime}
\subsetneq Z$, as required.
\bl
Conversely, suppose that there exists a proper subscheme $Z^{\prime}
\subsetneq Z$ such that the given extension defining $\F = \F_e$ is
induced from an extension $e^{\prime} \in  \text{Ext}^1(L \otimes
\I_{Z^\prime}, \O_X)$. Denote by $\E$ the double dual of the torsion-free sheaf
$\F_{e^\prime}$ determined by $e^\prime$, so that $\E$ is locally free of rank
two. Arguing from a diagram much like (3.11), one finds that there  exists
an injective sheaf homomorphism  $\nu : \F \lra \E$  with a non-trivial finite
cokernel.  But a generically injective map between two vector bundles of the
same rank is either an isomorphism or drops rank along a divisor (locally
defined as a determinant). Therefore $\F$ cannot be locally free.
\qed
\enddemo
\bl

We assume henceforth for simplicity that $Z$ is a {\it reduced} finite scheme,
i.e. that all the points of $Z$ appear with ``multiplicity one". Then we have:
\proclaim{Corollary 3.12}  There is an extension class $e \in
\text{\rm{Ext}}^1(L \otimes I_Z, \O_X)$ with $\F_e$ locally free if and only if
for every proper subset
$$Z^{\prime} \subsetneq Z \ \ (\text{including } Z^\prime = \emptyset),$$
$\text{\rm Ext}^1(L \otimes \I_{Z^\prime}, \O_X)$ maps to a proper subspace of
$\text{\rm Ext}^1(L \otimes \I_{Z}, \O_X)$.  \endproclaim
\demo{Proof} In fact, since $Z$ is reduced there are only finitely many proper
subschemes $Z^\prime \subset Z$. \qed \enddemo

\bl The next ingredient we'll need is:
\proclaim{Serre-Grothendieck Duality}  If $\G$ is any coherent sheaf on the
smooth projective surface $X$, then there is a natural isomorphism
$$
\text{\rm Ext}^1(\G, \O_X(K_X)) \lra H^1(X, \G)^* .$$ \endproclaim
\ni We refer to \cite{H2}, III.7, for the proof (cf. also \cite{GH1}, Chapter
5,
\S 4). However let us at least define the map appearing in the statement.
Represent a given element $e \in
\text{\rm Ext}^1(\G,
\O_X(K_X))$  by an extension $0 \lra \O_X(K_X) \lra \F_e \lra \G \lra 0$. Then
the connecting homomorphism defines a map
$$  H^1(X,\G) \lra H^2(X,\O _X (K_X))$$ which, after fixing an identification
$H^2(X,\O _X (K_X)) = \C$, may be viewed as an element of $H^1(X,\G)^*$. (See
also Exercise 3.19.(iv).)

\bl Still assuming that $Z$ is reduced, we may now state the answer to Question
3.7 as:
\proclaim{Theorem 3.13} {\cite{\rm \bf GH2}} There exists a rank two vector
bundle
$E$ with det$\ E = L$, plus a section $s \in \Gamma(E)$ with $Z(s) = Z$, if
and only if every section of
$\O_X(K_X + L)$ vanishing at all but one of the points of $Z$ also vanishes at
the remaining point.
\endproclaim
\bl \demo{Remark} When $Z = \{ x \}$ consists of a single (reduced) point, the
criterion is simply that $x$ be a base-point of $\O_X(K_X + L)$. \enddemo
\bl
\demo{Proof of Theorem 3.13}  As we observed above, the question is
equivalent to the existence of an extension (3.8) with $\F_e$ locally free.
Note to begin with the formal fact that it is enough to test the  condition in
Corollary 3.12 for sets
$Z^\prime
\subset Z$ obtained by deleting one point from
$Z$. Hence there exists a locally free extension of $L \otimes \I_Z$ by $\O_X$
if and only if for every point
$x \in Z$, the map
$$ \text{\rm Ext}^1(L \otimes \I_{Z-\{x\}} , \O_X) \lra \text{\rm Ext}^1(L
\otimes
\I_{Z} , \O_X) \tag * $$ is non-surjective. Now
$$
\gather
\text{\rm Ext}^1(L \otimes \I_{Z} , \O_X) = \text{\rm Ext}^1(\O_X(K_X + L)
\otimes
\I_{Z} , \O_X(K_X)), \\
\text{\rm Ext}^1(L \otimes \I_{Z-\{ x \}} , \O_X) = \text{\rm Ext}^1(\O_X(K_X +
L)
\otimes \I_{Z - \{ x \}} ,
\O_X(K_X)) .
\endgather
$$ So by Duality, the existence of a locally free extension is equivalent to
the
non-injectivity of
$$ H^1(X, \O_X(K_X + L) \otimes \I_Z) \lra H^1(X, \O_X(K_X + L) \otimes \I_{Z
- \{ x \}}) \tag **
$$ for every $x \in Z$. But the map (**) sits in the long exact sequence on
cohomology determined by the sheaf sequence
$$ 0 \lra \O_X(K_X + L) \otimes \I_Z \lra \O_X(K_X + L) \otimes \I_{Z  - \{ x
\}}
\lra
\O_X(K_X + L) \otimes \O_{ \{ x \}} \lra 0 .$$ In particular, the
non-injectivity of (**) is equivalent to evaluation at $x$ giving the zero
homomorphism $H^0(\O_X(K_X + L) \otimes \I_{Z  - \{ x \}}) \lra H^0(\O_X(K_X
+L) \otimes
\O_{ \{ x \}}) = \C$. But this means exactly that every section of $\O_X(K_X +
L)$ vanishing  on $Z - \{ x \}$ also vanishes at $x$, as claimed. \qed \enddemo
\bl

In the situation of Reider's Theorem 2.1(i), we have now achieved the goal of
encoding  geometric hypotheses on the linear series $| K_X + L |$ into the
existence of a (hopefully!) interesting  vector bundle on $X$:

\proclaim {Corollary 3.14} Let $L$ be any line bundle on the smoooth surface
$X$, and suppose that
$x \in X$ is a point at which every section of $\O_X(K_X + L)$ vanishes.  Then
there exists a rank two vector bundle $E$ on $X$, with det$ \ E = L$, which has
a section
$s \in \Gamma (X, E)$ vanshing precisely at $x$, i.e. $Z(s) = \{ x \}$. \qed
\endproclaim
\bl
\ni It remains to use this bundle to produce the required curve $D$. We take
this up in the next section. We conclude the present discussion with several
exercises, some of which outline other applications of vector bundles to study
linear series.

\bl \ni \xca{Exercise 3.15} ({\bf Non-Reduced Zero-Schemes}.) Generalize
Theorem
3.13 to  allow for the possibility that $Z$ is not reduced. In fact, assuming
that $Z \subset X$ is a finite local complete intersection subscheme, show that
the statement of Corollary 3.12 remains true if one deals with proper
subschemes $Z^\prime \subset Z$, with an analogous modification in the
statement
of Theorem 3.13.  Deduce that if $L$ is a  line bundle on $X$ such that
$\O_X(K_X + L)$ fails to be very ample, then there is a rank two vector bundle
$E$ with det$ \ E = L$ with a section vanishing on a finite scheme $Z \subset
X$ of length two. [The crucial point is that a local Gorenstein ring -- in
particular, a local ring of Z -- has a unique minimal non-zero ideal of
dimension one, to wit the socle (cf. \cite{C}). Hence a finite local complete
intersection scheme does not contain continuous families of maximal
proper subschemes.]
\qed
\endxca

 \ni \xca{Exercise 3.16} ({\bf Cayley-Bacharach Theorem}, \cite{GH2}.)  Let $X$
be a smooth surface, and let $C_1 , C_2 \subset X$ be effective reduced
divisors on $X$ which meet transversely. Prove that if $D \in |K_X + C_1 +
C_2|$ is a divisor passing through all but one of the points of $C_1
\cap C_2$, then $D$  passes through the remaining one as well. Generalize to
the
case when the intersection of
$C_1$ and $C_2$ is proper but possibly not transversal. \qed
\endxca

\xca{Exercise 3.17} ({\bf Elementary Transformations.})  Let $X$ be a smooth
surface, let $V$ be a vector bundle of rank $e$ on $X$, and let $C \subset X$
be a smooth curve. We suppose given a line bundle $A$ on $C$ of degree $d$,
plus a surjective map $\lambda : V|C \lra A$ from the restriction of $V$ to $C$
onto $A$.  From these data, we construct a new rank $e$ bundle $F$ on $X$, as
follows. We may view the invertible $\O_C$-module $A$ as a torsion $\O_X$ -
module; i.e. via ``extension by zero", $A$ becomes a coherent sheaf on $X$. We
denote this sheaf by
$\A$ to emphasize that it is not localy free on $X$. The composition
$V \lra V | C \lra A$ determines a surjection of $\O_X$-modules $\bar \lambda :
V
\lra \A$, and we set
$F = \ker \ \bar \lambda$. Thus one has the basic exact sequence:
$$ 0  \lra F \overset{\mu} \to  \lra V \lra \A \lra 0 \tag 3.17.1
$$ of sheaves on $X$.
\sbl
 (i).  Prove that $F$ is locally free on $X$, of rank $e$.  By analogy with a
classical construction on ruled surfaces (cf. \cite{H2}, V.5.7.1), $F$ is
called
the {\it elementary transformation} of $V$ determined by $\bar \lambda : V \lra
\A$.

\sbl
 (ii). Show that the Chern classes of $F$ are given by
$$
\aligned c_1(F) &= c_1(V) - [C] \\ c_2(F) &= c_2(V) - c_1(V) \cdot [C] + d,
\endaligned \tag 3.17.2
$$ where as above $d$ is the degree of $A$, considered as a line bundle on $C$.
[For
$c_1$ one can argue for instance that the map $\mu$ in (3.17.1) is a
homomorphism between two vector bundles of the same rank which drops rank on
$C$. The formula for
$c_2$ takes more work -- for example, one could deduce it by working
backwards from Riemann-Roch.]
\sbl
 (iii). Show that the transpose $\mu ^{*}$ of $\mu$ sits in the exact sequence
$$ 0 \lra V^{*} \overset{\mu^*} \to \lra F^* \lra \Cal{B} \lra 0, \tag 3.17.3$$
where $\Cal{B}$ is the extension by zero of the line bundle $B = N_{C/X}
\otimes A^* = \O_C(C)
\otimes A^*$ on $C$. [This amounts to the assertion that
${\Cal{E}xt}^1_{\O_X}(\A ,\O_X) = \Cal{B}$. When $A = \O_C$ this is elementary,
and in general one can argue that given a point
$P \in C$, the statement is true for $A$ if and only if it holds for $A(P)$.]
\qed
\endxca

\bl
\xca{Exercise 3.18} ({\bf Special Divisors on Curves on a K3 Surface.}) Let $X$
be a K3 surface, and let $C \subset X$ be a smooth curve of genus $g$. In this
exercise we outline a proof of the following
\sbl
\ni{ \bf Theorem} { {\cite{L1}}. {\it Assume that every curve in the linear
series {\rm $|C|$ } is reduced and irreducible. Then for every line bundle {\rm
$A$  } on  {\rm
$C$ }, one has  {\rm  $$g(C) - h^0(A)\cdot h^1(A) \ge 0.$$ } }

 To give some background for this statement, suppose for a moment that $C$ is
any smooth curve of genus $g$. It is important to understand under what
conditions on
$r, \ d,
\text{ and } g$ one can find a line bundle $A$ on $C$ of degree $d$ with
$h^0(A) \ge r+1$.  Classical parameter counts show that the family of all such
line bundles has expected dimension
$
\rho = \rho( d, g, r) =_{def} g - (r+1)(g - d +r)  \ge 0.
$ The existence of  $g^r_d$'s when $\rho \ge 0$ was established by Kempf and
Kleiman-Laksov in 1972. Classically, it was conjectured that if $\rho < 0$,
then a {\it general} curve $C$ of  genus $g$ carries no line bundles of degree
$d$ and $h^0
\ge r+1$. By Riemann-Roch, it is equivalent to predict:
$$\gathered
\text{ If } A \text{ is any line bundle on a general curve } C \text{ of genus
$g$, then } \\ g - h^0(A) \cdot h^1(A) \ge 0.
\endgathered
\tag *
 $$ Griffiths and Harris proved this (and more) in 1980 using a rather
elaborate
degenerational argument. We refer to \cite{Mfd} for a quick introduction to
this circle of ideas, and to \cite{ACGH} for a comprehensive overview and
references. Returning to the Theorem above, the hypothesis on the linear series
$| C |$ is certainly satisfied if $Pic(X) = \Z \cdot [C]$. On the other hand,
one knows from Hodge theory that for any integer
$g \ge 2$, there exists a K3 surface $X$ whose Picard group is generated by the
class of a curve
$C \subset X$ of genus $g$. Thus the Theorem gives a quick proof of (*),
without
degenerations.
\bl
  (i). Returning to the situation of the Theorem, let $A$ be a  line bundle of
degree $d$ on
$C$, with $h^0(A) =  r + 1$, such that both $A$ and $K_C \otimes A^*$ are
globally generated. Then $A$ is a quotient of the trivial vector bundle
$\O_C^{r+1}$ on $C$, so we can make an elementary transformation of $V =
\O_X^{r+1}$ to create a vector bundle $F = F(C,A)$ on $X$, of rank $r+1$,
sitting in the exact sequence
$0 \lra F \lra \O_X^{r+1} \lra \Cal{A} \lra 0.$ Set $E = F^*$. Prove that
$H^1(E)  = H^2(E) = 0$, and show that $E$ is generated by its global sections.
\sbl
 (ii). Prove that the holomorphic Euler characteristic of $E \otimes E^*$
satisfies
$\chi( E \otimes E^*) = 2h^0( E \otimes E^*) - h^1( E \otimes E^*) = 2 - 2
\rho(A),$ where $\rho(A) = g(C) - h^0(A) \cdot h^1(A)$. [Use the
Hirzebruch--Riemann--Roch formula
$
\chi( E \otimes E^*) = \int  Td(X) \cdot ch(E \otimes E^* )$, together with the
multiplicativity of the Chern character.]
\sbl
  (iii).  Assume now that $\rho(A) < 0$. Then $E$ has an endomorphism $w$ which
is not a multiple of the identity. Use this to construct a homomorphism $u : E
\lra E$ which drops rank everywhere on $X$. [If $\lambda$ is an eigenvalue of
$w(x)$ for some $x \in X$, put $u = w -
\lambda \cdot 1$.] Then consider the exact sequence
$0 \lra \text{im} \ u \lra E \lra \text{coker} \ u \lra 0$. Using the fact that
$E$ is globally generated, show that in the Chow group $A_1(X) = Pic(X)$, $c_1(
\text{im} \ u)$ and
$c_1(\text{coker} \ u)$ are represented by non-zero effective curves. Then
deduce  that $|C|$ contains a reducible or multiple curve.
\sbl
 (iv). The previous step proves the Theorem when both $A$ and $K_C \otimes A^*$
are globally generated. Show that the general case of the Theorem reduces to
this one. \qed
\endxca

\bl
\xca{Exercise 3.19} ({\bf Alternative Approach to Theorem 3.13}).  We  indicate
here a proof of the existence statement in Theorem 3.13 which avoids explicit
use of the sheaf-theoretic considerations appearing in the proof of Proposition
3.9.  As in the statement of the Theorem, let $Z \subset X$ be a reduced finite
scheme, say  $Z = \{ x_1, \dots , x_r \}$. Let
$$ f : Y = \text{Bl}_Z(X) \lra X$$ be the blowing-up of X along $Z$, and let $D
\subset Y$ be the exceptional divisor. Thus $D =
\sum_1^r D_i$, where $\P^1 \cong D_i \subset Y$ is the $(-1)$-curve lying over
$x_i
\in X$. The idea is that if the pair $(E,s)$ sought in (3.7) exists, then since
$f^*(s) \in \Gamma(Y, f^* E)$ vanishes on $D$, the pull-back of (3.6)
determines an extension $0 \lra \O_Y(D) \lra f^*E \lra f^*L(-D) \lra 0 \ $ on
$Y$. This suggests that we consider on $Y$ extensions of $f^*L(-D)$ by
$\O_Y(D)$. So fix $e \in \text{\rm Ext}^1(f^*L(-D), \O_Y(D))$ corresponding to
the exact sequence
$$ 0 \lra \O_Y(D) \lra  V \lra f^*L(-D) \lra 0. \tag *$$
\sbl (i).  Let $\P^1 \cong D_i \subset Y$ be one of the exceptional curves.
Show
that $V | D_i$ is an extension of $\O_{\P^1}(1)$ by $\O_{\P^1}(-1)$. Hence
either $V | D_i \cong \O_{\P^1}^2$, or $V | D_i \cong \O_{\P^1}(-1) \oplus
\O_{\P^1}(1)$.
\sbl (ii). Prove that if $V | D_i \cong \O_{\P^1}^2$ for every $i$, then $V =
f^*E$, where $E$ is a rank two vector bundle on $X$ with a section vanishing on
$Z$. So in this case we are done. [See
\cite{OSS}, I.2.2.6.]
\sbl (iii). Suppose  that $V |D_i \cong \O_{\P^1}(-1) \oplus \O_{\P^1}(1)$ for
some
$i$. Then make an elementary transformation along the resulting map $\bar
\lambda : V \lra \O_{\P^1}(-1)$ to define a new vector bundle $V^\prime$ on
$Y$. Using the fact that the composition $\O_Y(D) \lra V \lra
\O_{\P^1}(-1)$ is surjective, show that $V^\prime$ can be realized as an
extension of $f^*L(-D)$ by $\O_Y(D - D_i)$ from which (*) is induced. Thus the
given extension class $e$ lies in the image of $\text{\rm Ext} ^1(f^*L(-D),
\O_Y(D-D_i)) \lra
\text{\rm Ext} ^1(f^*L(-D), \O_Y(D)).$ (This is the analogue of Proposition
3.9.)
\sbl (iv). Recalling that $K_Y = f^* K_X + D$, reprove the existence
statement in Theorem 3.13 using Serre duality for line bundles on $Y$. \qed
\endxca

\bl
\bl
\head{\S 4. Reider's Theorem via Vector Bundles}\endhead

\bl
In the previous section we constructed a vector bundle $E$ encoding
the failure of an adjoint linear series $|K_X + L|$ to be very ample or
free. The next step is to study the geometry of $E$. Reider's basic tool
to this end is Bogomolov's instability theorem. We start with a quick
review of this fundamental result. Then we present the proof of Reider's
theorem.
\bl
\bl
%\ni {\it Bogomolov's Instability Theorem}.
\subhead{Bogomolov's Instability Theorem.}\endsubhead

\sbl Let $X$ be a smooth complex projective surface, and let $E$ be a
rank two vector bundle on $X$. It is always easy to construct very
negative rank one subsheaves of $X$: in fact, if $H$ is an ample
divisor, then for all
$n \gg 0$ there exist sheaf monomorphisms $\O_X(-nH) \hookrightarrow E$.
The notion of instability refers to the exceptional situation in which
$E$ has an unusually positive subsheaf. Bogomolov's theorem is a
numerical criterion for instability in terms of the Chern numbers of
$E$.

\bl We start with some formal definitions. Let $N(X)$ be the
N\'eron-Severi vector space of $X$, i.e. the subspace of $H^2(X, \R)$
generated by the classes of algebraic curves. The Hodge index theorem
implies that the intersection form on $N(X)$ has type $(+, - , \dots ,
-)$, so
$$ C = \{ \alpha \mid \alpha^2 > 0 \} \subset N(X)$$ is a cone with two
connected components. The {\it positive cone} $N(X)^+$ of $X$ is by
definition the component of $C$ containing the classes of ample divisors.

\demo{Definition 4.1} Let $E$ be a rank two bundle on $X$.  One says
that $E$ is {\it Bogomolov unstable} if there exist a finite subscheme
$Z \subset X$ (possibly empty), plus line bundles $A$ and $B$ on
$X$ sitting in an exact sequence
$$ 0 \lra A \lra E \lra B \otimes \I_Z \lra 0, \tag 4.1.1
$$  where  $$c_1(A) - c_1(B) \in N(X)^+. \tag 4.1.2 $$ \enddemo
\sbl
\ni Very concretely, (4.1.2) is equivalent to the conditions that $(A -
B)^2 > 0$ and that $(A - B)
\cdot H > 0$ for all ample divisors  $H$.  Note that the finite scheme
$Z$ in (4.1.1) is the scheme defined by the vanishing of the vector
bundle map $A \lra E$. It is suggestive to think of (4.1.2) as meaning
 roughly speaking that ``$A$ is more positive than
$B$".

\bl Bogomolov's statement is the following:
\proclaim {Theorem 4.2} {\rom{([{\bf Bog}] cf. also [{\bf Reid}].)}}
Let $E$ be a rank two vector bundle on a smooth projective surface
$X$. If $$ c_1(E)^2 - 4c_2(E) > 0, \tag 4.2.1 $$ then E is Bogomolov
unstable. \endproclaim
\bl
\ni We will outline a proof of Bogomolov's theorem from \cite{FdB1} in
\S6. We refer to
\cite{Bog}, \cite{Gies} or \cite{Mka} for other arguments, as well
as the corresponding statement for bundles of higher rank.
Shepherd-Barron \cite{S-B} has studied the geometry of this result in
positive characteristic, and Moriwaki \cite{Mwk} has given some
arithmetic analogues.

\xca{Exercise 4.3} If $E$ is a rank two bundle sitting in the exact
sequence (4.1.1), show that the Chern classes of $E$ are given by:
$$ c_1(E) = c_1(A) + c_1(B); \qquad c_2(E) = c_1(A) \cdot c_1(B) +
\text{ length}(Z).
$$ [For $c_2$, note that one can think of (4.1.1) as determining a map
$\O_X \lra E \otimes A^*$ with zero-scheme $Z$. In particular, $c_2(E
\otimes A^*) = \text{ length}(Z)$.] \qed
\endxca
\bl
\xca{Exercise 4.4}  Remark that if $E$ is a rank two vector bundle on
$X$, and $N$ is a line bundle, then $E$ is Bogomolov unstable if and
only if $E \otimes N$ is. Prove that up to scalar multiples,
$c_1^2 - 4c_2$ is the only weight two polynomial in  the Chern classes
of $E$ which is invariant under twisting by line bundles. Thus the
essential content of Theorem 4.2 is that there is {\it some} numerical
criterion for instability. \qed \endxca

\bl
\bl
%\ni {\it Proof of Reider's Theorem.}
\subhead{Proof of Reider's Theorem.} \endsubhead

\sbl We now turn to the proof of Reider's theorem. We follow the approach
of \cite{BFS}, which simplifies somewhat Reider's original presentation.
\footnote"*"{We remark that the argument given in \cite{L2} contains an
error.} We limit ourselves to the first statement of Reider's result,
leaving the second as an exercise for the reader.

\demo{Proof of Theorem 2.1.(i)} Fix a point  $x \in X$, and suppose that
$L$ is a nef line bundle on
$X$, with $L\cdot L \ge 5$, such that $\O_X(K_X + L)$ has a base-point
at $x$. By Corollary 3.14, there exists a rank two vector bundle $E$
sitting in the exact sequence
$$ 0 \lra \O_X \lra E \lra L \otimes \I_{x} \lra 0. \tag *
$$ Note from (4.3) that det$ \ E = L$ and $c_2(E) =  1$, and hence
$$c_1(E)^2 - 4c_2(E) = L\cdot L - 4 > 0.$$ Therefore Bogomolov's theorem
applies, and  one has an exact sequence
$$ 0 \lra A \lra E \lra B \otimes \I_Z \lra 0, \tag **
$$ where $Z$ is some finite subscheme of $X$, and $c_1(A) - c_1(B) \in
N(X)^+$.

\bl The plan now is to play off (*) against (**). Taking determinants,
we find in the first place that $A + B = L$, whence $A - B = 2A - L$.
Hence by definition of the positive cone:
$$
\aligned (2A - L)^2 &> 0 \\ (2A - L) \cdot H &> 0 \quad \forall \ \
\text{ample divisors } H.
\endaligned
\tag 4.5
$$ Denote by $$\alpha : A \lra L \otimes \I_{x}$$ the composition  $A
\hra E \lra L \otimes \I_x$ determined by (**) and (*). We claim to begin
with that $\alpha \ne 0$. In fact, in view of (*) it is enough to show
that Hom$(A, \O_X) = 0$, and this follows from the nefness of $L$ and the
second equation in (4.5). Then $\alpha$ is given by multiplication by
(the equation of) an effective divisor $D \subset X$, with
$$ x \in D  \qquad \text{and} \qquad A = L-D. $$
 It remains to show that $D$ satisfies the numerical conclusions of
Theorem 2.1.
\bl
To this end, we collect various inequalities. First:
$$ (L - 2D) \cdot L > 0. \tag 4.6$$ In fact, note that $L - 2D = 2A -
L$. Since $L$ is nef, and hence a limit of ample divisors, (4.5) implies
that in any event $(L - 2D) \cdot L \ge 0$. But if $(L - 2D) \cdot L =
0$, then $(L - 2D)^2 < 0$ by Hodge Index, and this contradicts (4.5).
\bl
 Next:
$$
\gather  (L^2)(D^2) \le (L \cdot D)^2 \tag 4.7 \\ (L-D)  \cdot D \le 1.
\tag 4.8
\endgather
$$ Indeed (4.7) is a consequence of Hodge Index. As for (4.8), computing
from (**) via (4.3), one finds $c_2(E) = (L-D)\cdot D + \text{ length} \
Z$. But $c_2(E) = 1$ and $\text{ length} \ Z \ge 0$, and this gives
(4.8). Finally, we claim:
$$ 2D^2 < L \cdot D. \tag 4.9 $$ Here we argue in cases. If $D^2 > 0$,
then $L \cdot D \ne 0$ by Hodge index. Moreover
$$ 2(L \cdot D)(D^2) \underset{(4.6)} \to{<} (L^2)(D^2) \underset{(4.7)}
\to{\le} (L \cdot D)^2,$$ and (4.9) follows. Next, say $D^2 = 0$. Then
Hodge index rules out the possibility that $L \cdot D = 0$.  Therefore
$L \cdot D > 0$, so again (4.9) is verified. Finally, if $D^2 < 0$, then
(4.9) is trivial since $L \cdot D \ge 0$.

\bl Combining (4.8) and (4.9), one finds:
$$ L \cdot D - 1 \le D^2 < \frac{L \cdot D}{2}.
$$ But this is only possible if
$$
\gather L \cdot D = 0 \ , \ D^2 = -1; \ \ \text{or} \\ L \cdot D = 1 \ ,
\ D^2 = 0.
\endgather
$$ This completes the proof. \qed \enddemo

\bl
\xca{Exercise 4.10} ({\bf Mumford's Proof of Vanishing.}) Let $X$ be a
smooth projective surface, and let $L$  be a nef line bundle on $X$ such
that $L^2 > 0$. Prove that then $H^1(X, \O_X(K_X + L)) = 0$. [If not,
there exists a non-split extension of $L$ by $\O_X$, to which one can
apply Bogomolov's theorem.  See \cite{Reid} for details.] \qed \endxca

\bl
\ni{\bf Exercise 4.11.} ({\bf Higher Order Embeddings.}) This exercise
is concerned with the following result of Beltrametti and Sommese:
\proclaim{Theorem} {\rm (\cite{BS})} Let $X$ be a smooth surface, let
$L$ a nef line bundle on $X$, and let $d$ be a positive integer such
that $L ^ 2 > 4d$.
 Then either the restriction
$$e_Z : \Gamma(X, \O_X(K_X + L)) \lra \Gamma(Z, \O_Z(K_X + L))  $$ is
surjective for every subscheme $Z \subset X$ of length $d$, or else
there exists an effective divisor $D \subset X$ such that
$$ L \cdot D - d \le D^2 < \frac{1}{2} (L \cdot D). \tag *
$$
\endproclaim
\bl (i). Assume that there exists a finite subscheme $Z \subset X$ of
length $d$ such that $e_Z$ fails to be surjective. By induction, one can
assume that $e_{Z^\prime}$ is surjective for every proper subscheme $
Z^\prime \subset Z$. Show that then there exists a rank two vector
bundle $E$ on $X$, with det$  \ E = L$, having a section $s \in
\Gamma(X, E)$ such that $Z(s) = Z$.
\sbl (ii). Arguing as in the proof of Reider's theorem, construct a
divisor $D$ satisfying (*).
\sbl (iii). Let $L$ be an ample line bundle on $X$ such that $L^2 > 4d$
and $L \cdot C \ge 2d$ for every irreducible curve $C \subset X$. Show
that then the restriction  $e_Z$ is surjective for every finite
subscheme $Z \subset X$ of length $\le d$. \qed

\bl
\xca{Exercise 4.12} ({\bf Gonality of Complete Intersection Curves.})
Reider-type methods can sometimes be used to study linear series on
subvarieties of codimension $\ge 2$ in an ambient manifold. To
illustrate the approach, we consider here a very concrete    question in
classical curve theory: what is the least degree required to express a
complete intersection curve $C \subset \P^r$ as a branched covering $C
\lra \P^1$ of the Riemann sphere? The answer is given in the following
\proclaim{Theorem} Let $C \subset \P^r$ be a smooth complete
intersection of hypersurfaces of degrees $2 \le a_1 \le a_2 \le \dots
\le a_{r-1}.$ Let $A$ be a base-point free line bundle on
$C$, of degree $d$, with  $h^0(C, A) \ge 2$. Then
$d \ge (a_1 - 1) \cdot a_2 \cdot \dots \cdot a_{r-1}.$ \endproclaim
\bl
\ni  The idea of the argument is this: let $S \supset C$ be a general
complete intersection surface of type $(a_2, \dots, a_{r-1})$. As in
Exercise 3.18, one can associate to $A$ a rank two vector bundle $F$ on
$S$. One finds that if $d < (a_1 - 1)  a_2 \cdots  a_{r-1}$, then $F$ is
Bogomolov unstable. It is easy to get a contradiction provided one knows
that the destabilizing subsheaf is of the form $\O_S(k)$, but this
doesn't seem to be guaranteed. To remedy this, instead of working on a
surface we work on a complete intersection threefold $X \supset C$,
whose Picard group is controlled by the Lefschetz theorems. Related
results, proved by more classical methods, appear in
\cite{CL} and \cite{Bas}, and the Theorem also connects with
some of the conjectures in \cite{EGH}. Paoletti \cite{Paol} has
extended the techniques of this exercise to deal with certain
non-complete intersection curves. He proves the striking result that
under suitable numerical hypotheses, the gonality of a space curve $C
\subset \P^3$ is governed by its {\it Seshadri constant}, which roughly
speaking measures how positive the hyperplane bundle $\O_{P^3}(1)$ is in
a neighborhood of
$C$ (see \S 5).
\sbl (i). Put $\gamma = a_3 a_4 \cdots a_{r-1}$, and let $X \supset C$
be a general complete intersection threefold of type $(a_3, \dots,
a_{r-1})$. [If $r=3$ take $X = \P^3$ and $\gamma = 1$.] Let $f : Y \lra
X$ be the blowing-up of $C$, and let ${ E} \subset Y$ be the exceptional
divisor, with $\pi : {  E} \lra C$ the natural map. Consider on ${  E}$
the globally generated line bundle $\A = \pi^* A$. Choosing a base-point
free pencil in $\Gamma({  E}, \A)$, we define in the usual way a rank
two vector bundle $\F$ on $Y$ via the sequence $0 \lra \F \lra
\O_Y^2 \lra \A \lra 0$. Compute the Chern classes of $\F$.
\sbl (ii). Denote by $H$ the pull-back to $Y$ of the hyperplane divisor
on $X$, and for $0 \le
\eps \in \Q$ consider the  $\Q$-divisor $D_{\eps} = (a_2 + \eps)H - {
E}$.
 Show that $D =_{\text{def}} D_0 = a_2 H - { E}$ is globally generated
and that $D_\eps$ is ample if $\eps > 0$. Now assume that $d < (a_1 -
1)a_2\gamma$. Prove that then for $0 < \eps \ll 1$:
$$
\left ( c_1(\F)^2 - 4c_2(\F) \right ) \cdot D_\eps = (a_1 -
\eps)a_1a_2\gamma -4d > 0. \tag *
$$

\sbl (iii). Fixing $\eps$ for which (*) holds, an extension by Miyaoka
\cite{Mka} of Bogomolov's instability theorem implies that there
exists a rank one subsheaf $\L \subset \F$ such that
$(2c_1(\L) - c_1(\F))\cdot D_\eps \cdot D > 0$. Show that one can assume
that
$\L$ is locally free, and that the vector bundle map $\L \lra \F$ drops
rank (if at all) on a codimension two subset $ Z \subset Y$. Prove that
$\L = \O_Y(-tH - \mu {  E}) $ for some integers $t, \mu \in \Z$. [Recall
that $Pic(X) = \Z$ thanks to the Lefschetz Hyperplane Theorem.]
\sbl (iv). Now let $S \in |a_2 H - {  E}| = |D|$ be a general divisor,
so that $S$ is isomorphic to a complete intersection surface of type
$(a_2, \dots , a_{r-1})$ through $C$. Setting $F = \F | S$, show that
$c_2(F) = d$, and that the restriction to $S$ of the subsheaf $\L \hra
\F$ gives rise to an exact sequence $0 \lra \O_S(-sH) \lra F \lra
\O_S((s-a_1)H) \otimes I_W \lra 0$, where $s = t + \mu a_1$, and $W
\subset S$ is some finite subscheme. Use  instability  to prove that $2s
< a_1$. Then estimate $c_2(F)$ to deduce that $a_1 < s+1$. But $s > 0$
since $h^0(S, F)  = 0$, and this gives a contradiction.
\sbl (v). Prove that the inequality in the Theorem is the best possible,
in the sense that for any integers $2 \le a_1 \le \dots \le a_{r-1}$,
there exists a complete intersection curve $C$ that carries a base-point
free pencil of degree $(a_1 - 1) \cdot a_2 \cdot \dots \cdot a_{r-1}$.
\qed
\endxca

\bl
\bl
\head{\S 5. Vanishing Theorems and Local Positivity}\endhead

\bl
The vector bundle methods described above seem essentially limited to
surfaces. In higher dimensions, vanishing theorems are the only tools that
have so far achieved significant success, and they will be the focus of
the rest of these lectures. In this section we discuss the questions that
arise if one tries to use vanishing in the most naive way to produce
pluricanonical or adjoint divisors. In a word, one is led to study the
``local positivity" of ample line bundles. We will see that one cannot
hope to recover completely the known results (e.g. Corollary 2.7) in this
fashion. However it turns out somewhat surprisingly that there are bounds
on local positivity which apply at a {\it generic} point of a smooth
surface, and these give results which in some respects go beyond Reider's
theorem. While we try to explain in  detail how vanishing theorems come
into the picture, we content ourselves with just a sketch of the
statements on local positivity.

\bl Let $X$ be a smooth projective surface, and let $L$ be an ample (or
nef and big) line bundle on $X$. We remarked in Proposition 1.2 that
$\O_X(K_X + L)$ is free at a point $x \in X$ if and only if
$H^1(X, \O_X(K_X + L)\otimes \I_x) = 0$. Standard vanishing theorems can't
directly apply here because the sheaf in question isn't locally free. The
traditional first step around this problem is to blow up at $x$, which at
least reduces the question to one involving only invertible sheaves.
\bl
	So fix a point $x \in X$,  let
$$f : Y =  Bl_{x}(X)  \lra X$$
be the
blowing up of $X$ at $x$, and denote by $E \subset Y$ the exceptional
divisor. The first point to note is:
\proclaim{Lemma 5.1} Let $L$ be a line bundle on the smooth surface $X$,
and let $r > 0$ be any positive integer. Then for all $i \ge 0$ there are
 isomorphisms
$$ H^i(X, \O_X(K_X + L) \otimes \I_x^r) = H^i(Y, \O_Y(K_Y + f^* L -
(r+1)E)).
$$
\endproclaim
\demo{Proof}  Observe first that  $f_*(\O_Y(-rE)) =  \I_x^r$, whereas
$R^jf_*(\O_Y(-rE)) = 0$ for $j > 0$. In fact, via the inclusion
$f_*(\O_Y(-rE)) \subseteq f_*(\O_Y) = \O_X$, we may identify the stalk of
$f_*(\O_Y(-rE))$ at $x$ as consisting of germs of functions on $X$ whose
pull-backs to $Y$ vanish to order $\ge r$ along the exceptional divisor,
which yields the first assertion. The second may be proven inductively by
taking direct images of the sequence $0 \lra \O_Y(-rE) \lra \O_Y(-(r-1)E)
\lra O_E(-(r-1)E) \lra 0.$ (Compare \cite{H2}, V.3.4.) Now recall that $K_Y
= f^*K_X + E$. Using the Leray spectral sequence and the projection
formula, one  finds:
$$
\align H^i(X, \O_X(K_X + L) \otimes \I_x^r) &= H^i(Y, f^* (O_X(K_X + L))
\otimes \O_Y(-rE)) \\
	&= H^i(Y, \O_Y(K_Y + f^* L - (r+1)E)),
\endalign
$$ as claimed. \qed \enddemo
\bl Recall next the statement of the basic:
\proclaim {Vanishing Theorem 5.2}  Let $V$ be a smooth projective variety.
If $L$ is any ample line bundle on $V$, then
$$H^i(V, \O_V(K_V + L)) = 0 \quad \text{for all } i > 0 .  $$
 More generally, the same statement holds assuming only that $L$ is nef
and big. \qed
\endproclaim
\bl
\ni The first assertion is of course the classical Kodaira vanishing
theorem. The fact that it is enough that $L$ be nef and big was proven by
Kawamata \cite{K1} and Viehweg \cite{V}. As we will see, it is very
convenient in applications only to have to check this weaker condition. We
refer to \cite{SS}, Chapter VII, or \cite{Kol1} for  nice introductions to
this extension  of Kodaira vanishing.

\bl The plan is to try to apply vanishing to  the line bundle $f^*L -
(r+1)E$ on $Y$. This requires knowing something about its positivity.
Demailly \cite{De2} has introduced a useful accounting mechanism for
keeping track of what one needs:
\definition{Definition 5.3} With notation as above, the {\it Seshadri
constant} of a nef line bundle $L$ at $x$ is the real number
$$\eps(L, x) = \sup \left \{ \eps \ge 0 \mid f^*L - \eps \cdot E \text{ is
nef} \right \}.$$
\enddefinition
\bl \ni Here $f^*L - \eps E$ is considered as an $\R$-divisor on $Y$, and
to say that it is nef means simply that $f^*L \cdot C^\prime \ge \eps( E
\cdot C^\prime)$ for all irreducible curves $C^\prime \subset Y$. Needless
to say, one can make the analogous definition on a smooth projective
variety of any dimension.

 \bl The connection with Seshadri's criterion for ampleness (cf.
\cite{H1}, Chapt. 1) occurs via:
\proclaim{Lemma 5.4} One has $$\eps(L,x) = \inf_{C \owns x} \left \{ \frac
{L \cdot C} { \text{ \rm mult}_x(C) }  \right \},$$ where the infimum is
taken over all reduced and irreducible curves $C
\subset X$ passing through $x$. \endproclaim
\demo{Proof} In fact, let $C \subset X$ be a reduced and irreducible curve
with multiplicity $m$ at
$x$, and let $C^\prime \subset Y$ denote the proper transform of $C$. Then
$C^\prime \cdot E = m$, and consequently for any $\eps$:
$$(f^*L - \eps E) \cdot C^\prime = L \cdot C - \eps m.$$ Hence if $ f^*L -
\eps E$ is nef, then $\eps \le \frac{L \cdot C}{m}$. The reverse
inequality is similar.
\qed \enddemo

\bl It is suggestive to think of the Seshadri constant $\eps(L, x)$ as
measuring how positive $L$ is locally near $x$. We present two exercises
that convey this point, and refer to \cite{De2}, \S 6,  for other
interpretations.
\bl
\xca{Exercise 5.5} (i). Show that if $L$ is a {\it very} ample line bundle
on the surface $X$, then $\eps(L, x)
\ge 1$ for all points $x \in X$.
\sbl (ii). Prove that the inequality in (i) holds assuming only that $L$
is ample and globally generated.
\sbl
\sbl (iii). Show that for any ample line bundle $L$ on $X$,  there exists a
positive constant $\eps = \eps(L) > 0$ such that  $\eps(L, x) \ge \eps$
for all $x
\in X$. (This is the elementary half of Seshadri's criterion for
ampleness.)
\qed
\endxca

\bl For the next exercise, and subsequent discussion, we need a
definition. Given a line bundle $B$ on $X$, and an integer $s \ge 0$, we
say that the linear series $|B|$ {\it generates}  $s${\it -jets} at $x$ if
the natural map
$$H^0(X, B)  \lra H^0(X, B \otimes \O_X/ \I_x^{s+1})$$ is surjective. In
other words, we ask that we be able to find a global section of $B$ with
arbitrarily prescribed $s$-jet at $x$. For instance $|B|$ generates
$0$-jets at every point of $X$ if and only if it is free, and $|B|$
generates  all $1$-jets if and only if the differential
$d\phi_{|B|}$ of the corresponding map $\phi_{|B|} : X \lra \P$ is
everywhere injective.
\bl
\xca{Exercise 5.6} (\cite{De2}, Theorem 6.4). Given a  line bundle $B$ on
$X$, let  $s(B , x)$ be the largest integer such that the global sections
in $H^0(X,B)$ generate
$s$-jets  at  $x \in X$. Returning to an   ample line bundle $L$
on $X$, put
$$ \sigma(L, x) = \limsup_{k \to \infty} \frac{1}{k} s(kL, x). $$ Show
that $\eps(L, x) = \sigma(L, x)$ for every $x \in X$. [Compare Proposition
5.10 below.]  \qed \endxca
\bl

We now have:
\proclaim{Proposition 5.7}  Let $L$ be an ample line bundle on the smooth
surface $X$, and let $x
\in X$ be a fixed point. If $\eps(L,x) > 2$, then $\O_X(K_X + L)$ is free
at $x$. More generally, if
$\eps(L,x) > s+2$ for some  integer $s \ge 0$, then the linear series
$|K_X + L|$ generates $s$-jets at $x$. The same statement holds if
$\eps(L,x) = s + 2$ provided that
$L^2 > (s+2)^2$. \endproclaim
\demo{Proof} To prove that $|K_X + L|$ generates $s$-jets at $x$, it is
sufficient to establish the vanishing
$$H^1(X, \O_X(K_X + L) \otimes \I_x^{s+1}) = 0.  \tag *$$  As before, let
$f : Y \lra X$ be the blowing up of $X$ at $x$, and let $E \subset Y$ be
the exceptional divisor. Setting $\eps = \eps(L,x)$, one has on $Y$ the
numerical equivalence of
$\R$-divisor classes:
$$f^*L - (s+2)E \equiv \frac{s+2}{\eps} (f^*L - \eps E) + (1 -
\frac{s+2}{\eps})f^*L.$$ The first term on the right is nef, and the
second is nef and big since $\eps > s+2$. Hence the vanishing theorem
(5.2) applies to $f^*L - (s+2)E$, and (*) then follows from Lemma 5.1. If
$\eps = s+2$, then $f^*L - (s+2)E $ is nef, and the inequality in the
statement of the Proposition implies that it is big. \qed \enddemo
\proclaim{Corollary 5.8} Let $A$ be an ample line bundle on the surface
$X$ such that $\eps(A,x) \ge 1$ at some point $x \in X$. Then for all $s
\ge 0$ the linear series $|K_X + (s+3)A|$ generates
$s$-jets at $x$. \qed \endproclaim
\bl
\xca{Exercise 5.9}  Show that if $\eps(L, x) > 4$ for all $x \in X$, then
$\O_X(K_X+L)$ is very ample. State and prove the analogues of (5.1), (5.7)
and (5.8) on a smooth projective variety $V$ of arbitrary dimension. \qed
\endxca

\bl As one might expect, it can happen that $\O_X(K_X + L)$ is free at $x$
without this being accounted for by the bound $\eps(L, x) \ge 2$. (See
Proposition 5.12  below.) However by a small variant of \cite{De2},
Theorem 6.4, statements such as (5.8) for jets of arbitrarily high order
are {\it equivalent} to inequalities on the Seshadri constant:
\proclaim{Proposition 5.10}  Let $A$ be an ample line bundle on $X$, $x
\in X$ a fixed point, and
$\eps > 0$ a real number. Suppose that for all $s \gg 0$, the linear
series $|K_X + rA|$ generates
$s$-jets as soon as $r > (s+2)/{\eps}$. Then $\eps(A, x) \ge \eps$.
\endproclaim
\demo{Proof} Let $C \owns x$ be a reduced and irreducible curve, with
mult$_x(C) = m$. Fix $s \gg 0$, and let $r = r(s)$ be the least integer $>
({s+2})/{\eps}$. Then $|K_X + rA|$ generates $s$-jets, and consequently we
can find a curve $D = D_s \in |K_X + rA|$, with mult$_x(D) = s$, having a
prescribed tangent cone at $x$. In particular, we can choose $D$ so that
the tangent cones to $C$ and $D$ at $x$ meet properly in $T_xX$, and since
$C$ is irreducible it follows that $C$ and $D$ themselves meet properly.
Therefore
$$C \cdot (K_X + rA) \ge \text{ mult}_x(C) \cdot \text{ mult}_x(D) = ms,$$
whence
$$\frac{C \cdot A}{m} \ge \frac{s}{r} - \frac{C \cdot K_X}{rm}.$$ But we
can assume that $r \le (\frac{s+2}{\eps}) + 1 $, and the result follows
 from Lemma 5.4 upon letting $s \to \infty$. \qed \enddemo

\bl So far this discussion has been quite formal. It remains to say
something about the actual behavior of the Seshadri constants $\eps(L, x)$.
As measures of local positivity, these are in any event very interesting
invariants, quite apart from the potential application to adjoint series.
We may summarize the story on surfaces in the following two statements:
\proclaim{Theorem 5.11} {\rm \cite{EL3}.} Let $L$ be an ample line
bundle on a smooth projective surface $X$. Then $\eps(L, x) \ge 1$ for all
except perhaps  countably many points $x \in X$, and moreover if
$L^2 > 1$ then the set of exceptional points is finite. If $L^2 \ge 5$ and
$L \cdot C \ge 2$ for all curves $C \subset X$, then $\eps(L, x) \ge 2$
for all but finitely many $x \in X$. \endproclaim

\proclaim{Proposition 5.12} {\rm ({\bf Miranda}).} Given $\eps > 0$, there
exists a surface $X$, a point $x \in X$, and an ample line bundle $L$ on
$X$ such that $\eps(L, x) \le \eps$. \endproclaim

\bl \ni It follows for instance from (5.7) and the  Theorem that if $L^2
\ge 5$ and $L \cdot C \ge 2$ for all $C \subset X$ then $\O_X(K_X + L)$ is
free off a finite set. But of course we know from Reider's theorem
(Corollary 2.6) that in fact $\O_X(K_X + L)$ is  {\it everywhere} globally
generated. So from this point of view, one may think of Seshadri constants
as giving local Reider-type results, which however are valid only at a
general point of the surface $X$. On the other hand, the statements coming
from (5.7) and (5.11) for higher order jets are necessarily stronger than
the uniform results deduced from Reider's method. (See Exercise 5.14 for a
summary.)

\bl We start by outlining Miranda's construction  of examples of small
Seshadri constants.
\demo{Proof of Proposition 5.12} Let $D \subset \P^2$ be an irreducible
plane curve of degree $d \gg 0$ with a point $x \in D$ of multiplicity
$m$. Let $D^\prime$ be a second irreducible curve of degree $d$, meeting
$D$ transversely. Choosing $D^\prime$ generally, we may suppose that all
the curves in the pencil spanned by $D$ and
$D^\prime$ are irreducible. Blow up the base-points of this pencil to
obtain a surface
$X$, admitting a map $\pi : X \lra \P^1$ with irreducible fibres, among
them $D
\subset X$. Observe that $\pi$ has a section $S \subset X$ (viz. an
exceptional curve of the blowing up $X \lra \P^2$) meeting $D$
transversely at one point. Fix an integer $a \ge 2$. It follows from the
Nakai criterion that the divisor  $L = aD + S$ on $X$ is ample. But $L
\cdot D = 1$ whereas mult$_x(D) = m$, so $\eps(L, x) \le \frac{1}{m}$.
Note that by taking suitable $a$ we can even make $L^2$ arbitrary large,
and by taking $L$ to be a multiple of $aD + S$ we can arrange that the
intersection numbers $L \cdot C$ of $L$ with irreducible curves $C \subset
X$ be bounded below by any preassigned integer. \qed \enddemo
\bl Finally, we sketch without full details the main idea of the proof of
Theorem 5.11. The argument is very elementary, but is essentially limited
to surfaces.
\demo{Idea of proof of Theorem 5.11} We focus on the first statement, and
we use the characterization (5.4) of Seshadri constants. The main point,
which was inspired by \cite{Xu}, is to view the question variationally. In
fact, the set
$$
\left \{ (C,x)  \mid  C \subseteq X
                \text{ is a reduced, irreducible curve with mult}_x(C)> L
\cdot C \right \}$$ consists of at most countably many algebraic families.
The first statement of the  theorem will follow if we show that each of
these families is discrete, i.e. that pairs $(C,x)$ forcing $\eps (L,x)<1$
are rigid.
\sbl Suppose to the contrary that $(C_t,x_t)$ is a non-trivial
one-parameter family of reduced and irreducible curves $C_t\subseteq X$
and points $x_t\in C_t$ with  mult$_{x_t}(C_t)>L\cdot C_t$. Let
$C=C_{t^*}$ and $x=x_{t^*}$ for general $t^*$, and set
$m=\text{mult}_{x_{t^*}}(C_{t^*})$. A  local computation involving
deformation theory of singular curves shows that
$$C^2\geq m(m-1).$$ (In brief, the given deformation of
$C$ determines a section $\rho(\frac{d}{dt}) \in H^0(\O_C(C))$, and one
shows that since the deformation preserves the $m$-fold point of $C$,
$\rho(\frac{d}{dt})$ vanishes to order $\ge m-1$ at $x$.)  But
$L^2\cdot C^2 \leq (L\cdot C)^2$ by the Hodge index theorem, and since
$L\cdot C\le  m-1$ by assumption, we find that
$$m(m-1)\leq L^2\cdot C^2 \leq ( C)^2 \leq (m-1)^2 .$$ But this is a
contradiction when $m>1$, and the first statement follows. \qed
\enddemo

\demo{Remark 5.13} Theorem 5.11 has recently been extended to varieties of
arbitrary dimension, although the numerical bound obtained is somewhat
weaker. In fact, it is shown in
\cite{EKL} that if $L$ is an ample line bundle on a smooth projective
variety of dimension $n$, then
$\eps(L, x) \ge \frac{1}{n}$ for all $x \in X$ outside a countable union
of proper subvarieties. The idea is to choose a divisor $D = D_x \in |kL|$
($k \gg 0$) with large multiplicity at $x$, and to study in effect the
higher order deformation theory of the pair $(D_x, C_x)$, where $C_x
\subset X$ is a Seshadri-exceptional curve at $x$. The argument is
inspired on the one hand by some techniques that come up in the theory of
diophantine approximation, and on the other hand by the methods used to
prove boundedness of Fano varieties of  given dimension. \enddemo
\bl
\xca{Exercise 5.14} ({\bf Higher Jets and Adjoint Series.}) It is
interesting to summarize the statements on separation of jets that come
out of these discussions. Let $A$ be an ample line bundle on the smooth
projectve  surface $X$. Prove the following:
\sbl (i).  The adjoint series $| K_X + (s+3)A |$ generates $s$-jets at a
sufficiently general point
$x \in X$.
\sbl (ii). There cannot exist a linear function $f(s)$ such that for all
$X$, $A$ and $s \gg 0$,
$|K_X + f(s)A|$ generates $s$-jets at {\it every} point $x \in X$.
\sbl  (iii). There exists a {\it quadratic} function $f(s)$ (independent
of $X$ and $A$) such that for $s
\gg 0$,  $|K_X + f(s)A|$ generates $s$-jets at every point $x \in X$. In
fact, for $s \ge 1$ one can take $f(s) = (s+1)(s+2)$. [Use Exercise
4.11.(iii).] \qed
\endxca

\bl
\xca{Exercise 5.15} ({\bf Normal Generation of ``Hyper-adjoint" series,}
\cite{BEL},\S3,
\cite{ABS}.) Let $V$ be a smooth projective variety of dimension $n$, and
let $B$ be a {\it very ample} line bundle on $V$.
\sbl
(i).  Prove that $K_V + (n+1)B$ is globally generated, and that $K_V +
(n+2)B$ is very ample. [For the first statement, consider a general divisor
$W \in |B|$ and argue by induction  on $n$.]
\sbl
(ii). Show that $\O_V(K_V + (n+1)B)$ is normally generated, i.e. that the
homomorphisms
$$Sym^m(H^0( \O_V(K_V + (n+1)B))) \lra H^0(\O_V(m(K_V + (n+1)B)))$$
are surjective for $m \ge 0$. [The case $m = 2$ is the essential one. Put
$N = \O_V(K_V + (n+1)B)$, and consider on $V \times V$ the line bundle $F =
pr_1^*(N) \otimes pr_2^*(N)$. It suffices to prove that $H^1(V \times V, F
\otimes \I_{\Delta} ) = 0$, where $\Delta \subset V
\times V $ is the diagonal. Using the hypothesis that
$B$ is very ample, show that $pr_1^*(B) \otimes pr_2^*(B) \otimes
\I_{\Delta}$ is globally generated. Now apply vanishing on the blow-up of
$V \times V$ along $\Delta$.] See \cite{EL1} for some generalizations
involving defining equations and higher syzygies, and \cite{BEL} for some
similar elementary applications of vanishing theorems to study the
equations defining projective varieties. \qed
\endxca

\bl
\xca{Exercise 5.16} ({\bf Seshadri Constants Along Finite Sets}.) This
exercise gives an outline of some unpublished work of  {\bf Geng  Xu}, and
we thank him for his permission to include it.  Let $X$ be a smooth
surface, and let $Z \subset X$ be a finite set consisting of $r$ points,
say $Z = \{ x_1, \dots , x_r \}$. Given an effective divisor $D \subset
X$, define mult$_Z(D) = \sum \ \text{mult}_{x_i}(D)$. If $L$ is a nef line
bundle on $X$, we define the Seshadri constant of $L$   along $Z$ to be
the real number
$$\eps(L,Z) = \inf_{C \subset X} \left \{ \frac
{L \cdot C} { \text{ \rm mult}_Z(C) }  \right \},$$
the infimum being taken over all reduced and irreducible curves $C
\subset X$ (not necessarily passing through $Z$).
\sbl
(i). Assume that $L$ is nef and $L^2 > r$. Show that if $Z$ consists of
$r$ sufficiently general points, then $\eps(L,Z) \ge 1$. [If not, let $C
\subset X$ be a Seshadri exceptional curve  which deforms with all
local deformations of the $x_i$.  Put $m_i = \text{mult}_{x_i}(C) \ge
0$. We may assume that $C$ passes through $x_1$, so $m_1 \ge 1$, and by
hypothesis $L \cdot C \le (m_1 - 1) + m_2 + \dots + m_r$. As $C$ moves
in an algebraic family generically covering $X$, $L \cdot C \ge 1$.
Consider the deformation  obtained by letting $x_1$ move and fixing $x_2,
\dots, x_r$. Since $C$ is reduced and irreducible, one finds as in
\cite{EL3},
\S1, the inequality:
$$C ^2 \ge (m_1 - 1)m_1 + m_2^2 + \dots + m_r^2 \ge (m_1 - 1)^2 + m_2^2 +
\dots + m_r^2.$$
Hodge index and the hypothesis $L^2 \ge r+1 $ then gives
$$(r+1) \left ( (m_1 - 1)^2 + m_2^2 + \dots + m_r^2 \right ) \le (m_1 - 1 +
m_2 + \dots + m_r)^2,$$
and by minimizing the the difference of the two sides one sees that this is
impossible.]
\sbl
(ii). In the situation of (i), suppose that $D \subset X$ is any
effective divisor (possibly non-reduced or reducible). Show that then
$\text{mult}_Z(D) \le L \cdot D$.
\sbl
\ni See Exercise 7.8  for an interesting application of this result. \qed
\endxca

\bl
\bl
\head{\S 6.  Adjoint  Series  and Bogomolov Instability via
Vanishing.}
\endhead

\bl We saw in the previous section that one can't hope to prove known
and expected results on  linear series using only the most naive
application of vanishing theorems for line bundles. Our purpose here
is to show how more subtle vanishing theorems, for $\Q$-divisors, do
lead to (parts of) Reider's theorem. This argument appears in
\cite{EL2}, \S1, and follows the approach pioneered by Kawamata, Reid,
Shokurov {\it et al.} in connection with the minimal model program.
Many of the complexities of the general KRS machine disappear on
surfaces, and the ideas become particularly transparent. We hope that
the present discussion can serve as a low-key introduction to this
important and powerful tool.  Completing this circle of
ideas, Fern\'andez del Busto \cite{FdB1} has shown that one can use
the approach of \cite{EL2} to give a new proof of Bogomolov's
instability theorem. We sketch the argument at the  end of the
section, and in Exercises 6.23 and 6.24.

\bl By way of motivation, let $L$ be an ample line bundle on a
smooth projective surface $X$, and consider the problem of constructing
a section of $\O_X(K_X + L)$ which is non-vanishing at some point $x
\in X$. The ``classical" approach might be to prove something along the
following lines:
\sbl
\xca{Exercise 6.1} Assume that $L^2 \ge 5$, and suppose there exists a
reduced irreducible divisor
$$D \in |L| \ \text{ such that } q = \text{ mult}_x(D) \ge 2.$$ Then
$H^1(X, \O_X(K_X + L) \otimes \I_x) = 0$, and consequently $\O_X(K_X +
L)$ is free at $x$. [Let
$f : Y \lra X$ be the blowing up of $X$ at $x$, with $E \subset X$ the
exceptional divisor. Let
$D^\prime \subset Y$ be the proper transform of $D$, so that $f^* D -
2E \equiv D^\prime + (q-2)E$ is effective and has positive
self-intersection. Prove that $f^* D - 2E$ is numerically
$1$-connected, and deduce from Ramanujam vanishing (\cite{BPV},
IV.8.2) that $H^1(Y, \O_Y(K_Y + f^*L - 2E)) = 0$. Then use Lemma 5.1
to conclude.] \qed \endxca
\bl
\ni The essential drawback to this statement is that in the situations
of interest there is no reason even to suppose that $|L|$ is
non-empty. The advantage of the KRS method is that it lets one make an
{\it asymptotic} construction. Specifically, we will take $k \gg 0$
and work with a divisor
$$ D \in |kL| \  \text{ such that mult}_x(D) > 2k.$$ There is no
problem in producing $D$ if $L^2 \ge 5$ and $k$ is sufficiently large.
Then one will want to ``divide $D$ by $k$", and this is where
$\Q$-divisors come into the picture.
\sbl
\bl
%\ni {\it Vanishing Theorems for $\Q$-Divisors }
\subhead{Vanishing Theorems for $\Q$-Divisors.}\endsubhead

\sbl We start with some notation and definitions. Let $V$ be a smooth
projective variety of dimension
$n$. A $\Q$-{\it divisor} on $V$ is simply a $\Q$-linear combination
of prime divisors:
$$ M = \sum a_i D_i \ \ (a_i \in \Q).$$ The multiplicity  mult$_x M$
of $M$ at a point $x \in X$ is taken to be $\sum a_i \text{mult}_x
D_i$. Assuming that the $D_i$ are distinct, we define the {\it
round-up} of $M$ to be the integral divisor $\rndup{M} = \sum
\rndup{a_i} D_i$, where $\rndup{a_i}$ denotes the least integer $ \ge
a_i$. The {\it integer part}, or {\it round-down}  $[M]$ is defined
similarly, and the {\it fractional part} $\{ M
\}$ of $M$ is $\{ M \} = M - [M]$. There is a $\Q$-valued
intersection theory involving $\Q$-divisors, defined in the evident
way by first clearing denominators, and one has the usual functorial
operations such as pull-backs under morphisms. This gives rise to
the notion of numerical equivalence of $\Q$-divisors, which we
continue to denote by $\equiv$. We say that $M$ is {\it nef} if $M
\cdot C \ge 0$ for all irreducible curves
$C \subset V$, and $M$ is {\it ample} if the statement of Nakai's
criterion holds.  Equivalently, $M$ is nef or ample if $mM$ is so,
for some positive integer $m>0$ such that $mM$ is an integral
divisor. If $M$ is nef, it is in addition {\it big} if
$(M^n) > 0$.

\bl The basic result is:
\proclaim{Kawamata-Viehweg Vanishing Theorem 6.2} Let $M$ be a nef and
big $\Q$-divisor on the smooth projective variety $V$. Assume that the
fractional part $\{ M \}$ of $M$ is supported on a divisor with global
normal crossings. Then
$$ H^i(V, \O_V(K_V + \rndup{M})) = 0 \ \text{ for } \ i > 0. \
$$ \endproclaim
\sbl
\ni In other words, the Theorem gives a vanishing for (integer)
divisors of the form:
$$   K_V  +  \left ( \text{ nef and big $\Q$-divisor } \right ) +
\Delta, $$ where $\Delta$ is an effective fractional divisor (i.e. $[
\Delta] = 0$) with normal crossing support. The statement may appear
non-intuitive at first blush, but we will see momentarily that this is
precisely the tool needed to generalize the argument sketched in
Exercise 6.1.  We won't prove (6.2) in these notes. The original
arguments of \cite{K1} and \cite{V} used covering constructions to
deduce the result from vanishing for integer divisors. The normal
crossing hypothesis is used to contol the singularities introduced
upon passing to a covering. A number of direct proofs have since been
given, one based on connections with logarithmic singularities
\cite{EV1}, another on Hodge theory for twisted
coefficient systems \cite{Kol1}, and a third involving singular metrics
on line bundles \cite{De2}.  We refer to  \cite{Kol1}, \cite{CKM},
Chapter 8,  \cite{Kol3}, Chapters 9 - 10, and \cite{EV2} for good
discussions. It is worth emphasizing that Theorem 6.2 is by now not
much harder to prove than the classical Kodaira vanishing theorem.
\bl
\demo{Remark 6.3}  Since rounding does not in general respect linear
equivalence, it is essential that the fractional part of the
$\Q$-divisor $M$ appearing in (6.2) be defined as an actual divisor,
and not merely as an element in $Pic(X) \otimes \Q$. However we will
often identify two $\Q$-divisors if their fractional parts coincide
and their integer parts are linearly equivalent. By the same token, we
will deal with ``hybrid" objects of the form $L+D$ where $L$ is a
line bundle (defined up to isomorphism) and $D$ is a $\Q$-divisor.
\enddemo
\bl

In the hope of conveying right away some feeling for how Theorem 6.2
is used, we prove a criterion extending Exercise 6.1. It asserts
roughly that the existence of a divisor in $|kL|$ with an ``almost
isolated" singularity of high multiplicity gives rise to a
non-vanishing section of
$\O_X(K_X + L)$.  This result will play an important role in our proof
of Reider's theorem. In the argument below, we will make use of the
fact --- to be established in Exercise 6.6
 --- that on a surface, one can ignore  the normal crossing hypothesis
in (6.2).
\proclaim{Proposition 6.4}  Let $L$ be a nef and big  line bundle on a
smooth projective surface $X$. Fix a point $x \in X$, and an integer
$s \ge 0$. Suppose that for some $k > 0$ there exists a divisor
$$D \in |kL| \ \text{  with $\ q =_{\text{def}}$ mult}_x(D) >
(s+2)k,$$ plus an open neighborhood $U \owns x$ of $x$ in $X$ such
that
$$\text{mult}_y(D) < q/(s+2) \ \text{ for all } \ y \in U - \{ x \} .
$$ Then
$$H^1(X, \O_X(K_X + L) \otimes \I_x^{s+1}) = 0,$$ i.e. $|K_X + L|$
generates $s$-jets at $x$.
\endproclaim

\demo{Proof}  Write $D = \sum d_i D_i$, where the $D_i$ are distinct
prime divisors. The upper bound on mult$_y(D)$ implies that if $D_i$
passes through $x$, then $d_i < q/(s+2)$. For simplicity we  assume
for the time being that every component of $D$ passes through $x$; the
changes necessary in general will be sketched at the end of the proof.
Consider as before the blowing-up $f : Y \lra X$ of $X$ at $x$. Let
$E \subset Y$ be the exceptional divisor, and denote by
$D_i^\prime \subset Y$ the proper transform of $D_i$. Thus $f^*D = qE
+ \sum d_i D_i^\prime$. The idea of the proof is to study the
$\Q$-divisor
$$M = f^*L - \left ( \frac{s+2}{q} \right ) f^* D = f^* L - (s+2)E -
\sum \left ( \frac{s+2}{q} \right ) d_i D_i^\prime.
\tag *
$$ The first point to note is that by hypothesis
$$\left ( \frac{s+2}{q} \right ) d_i< 1 \ \text{ for all }i,$$  and
therefore  $$K_Y + \rndup{M} = K_Y + f^*L - (s+2)E.$$
 On the other hand, one has the numerical equivalence:
$$M \equiv f^*L - \left ( \frac{s+2}{q}\right )f^* D \equiv \left (1 -
\left (\frac{s+2}{q}\right ) k
\right ) f^*L,$$ and as
$q > k(s+2)$, it follows that $M$ is nef and big. Since we are on a
surface we don't need to worry about normal crossings, and therefore
Theorem 6.2 gives $H^1(Y, \O_Y(K_Y + f^*L - (s+2)E)) = 0$. We conclude
with Lemma 5.1.
\sbl It remains to treat the possibility that not all of the
components of $D$ pass through $x$. In this case, we only know that
$d_i < q/(s+2)$ for those $i$ such that $D_i \owns x$. Defining $M$ as
in (*), it follows that $K_Y + \rndup{M} = f^*L - (s+2)E - N^\prime$,
where $N^\prime \subset Y$ is an effective (or zero) integral divisor
whose support is disjoint from $E$. Let $N = f_* N^\prime$ be the
corresponding divisor on $X$, so that $x \notin \text{supp} \ N$.
Then as above vanishing gives
$H^1(X , \O_X(K_X + L - N) \otimes \I_x^{s+1}) = 0$, i.e. we can find
a section of $\O_X(K_X + L -N)$ with arbitrarily prescribed $s$-jet at
$x$. But one has an inclusion $H^0(\O_X(K_X + L -N)) \subset
H^0(\O_X(K_X + L))$, and since $N$ doesn't pass through $x$, it
follows that $|\O_X(K_X + L)|$ also generates $s$-jets at $x$, as
required.
\qed
\enddemo

\bl
\demo{Remark 6.5}  There is an analogous statement in higher
dimensions:
\proclaim{\rm (6.5.1)} Let $V$ be a smooth projective variety of
dimension $n$, and let $L$ be an ample line bundle on $V$. Given a
point $x \in V$, assume that there exists a divisor $D \in |kL|$ with
$$\text{mult}_x D > (n+s)k$$ but mult$_yD \le k$ for $y$ in a punctured
neighborhood of $x$. Then $|K_X + L|$ generates $s$-jets at $x$.
\endproclaim
\bl
\ni We learned of (6.5.1) from Siu. Esnault and Viehweg give a proof
in \cite{EV2}, (7.5), (7.7). The argument is more involved than in the
surface case, because in general one has to pass to an embedded
resolution of $D$ in order to apply Kawamata-Viehweg vanishing. Note
also that in order to verify the upper bound on mult$_y(D)$  in
dimension three or higher, it is not enough to know simply that the
components of $D$ appear with low multiplicity. \enddemo
\bl
As we have already remarked, it is a happy fact that one can ignore
the normal crossing hypothesis in  Theorem 6.2 when working on
surfaces:
\xca{Exercise 6.6} ({\bf Sakai's Lemma}, \cite{Sak1}, cf. \cite{EL2},
\S1.) Let $X$ be a smooth surface, and let $M$ be any big and nef
$\Q$-divisor on $X$. Then
$$H^i(X, \O_X(K_X + \rndup{M})) = 0 \ \text{ for } i > 0.$$  [A simple
if somewhat inefficient argument proceeds as follows. By a sucession
of blowings up at points, one constructs a map $\phi : X_1 \lra X$
such that the fractional part $\{ \phi^* M \}$ of
$\phi^* M$ is supported on a normal crossing divisor, and hence
$H^i(X_1, \O_{X_1}(K_{X_1} + \rndup{ \phi^* M }) ) = 0$ for $i > 0$.
Working step by step down from
$X_1$, it is then enough to prove the following:
\ext  Let $f : Y \lra X$ be the blowing up of a smooth surface at a
point $x \in X$, and let $M$ be a
$\Q$-divisor on $X$. If  $H^i(Y, \O_Y(K_Y + \rndup{f^* M})) =  0$ for
some $i > 0$, then $H^i(X,
\O_X(K_X + \rndup{  M})) =  0$. \endext
\sbl
\ni In fact, show that
$$K_Y + \rndup{f^* M} =  f^*(K_X + \rndup{M}) - pE$$  for some $p \ge
-1$. Then
$$ H^i(Y , \O_Y(K_Y + \rndup{f^*M} )) =  H^i(X , \O_X(K_X + \rndup{M})
\otimes {\I}_x^p )
$$  for all $i$, where we make the convention that $\I_x^p = \O_X$ if
$p = -1$. The assertion follows.]
\qed
\endxca

\bl
\xca{Exercise 6.7} ({\bf Algebro-Geometric Multiplier Ideals},
\cite{EV2}, (7.4), (7.5).) Let $M$ be any nef and big $\Q$-divisor on
a smooth projective variety $V$ of dimension $n$. Prove that there
exists an ideal sheaf ${\Cal J}_M \subset \O_V$, with $\O_V/{\Cal
J}_M$ supported in codimension $\ge 2$, such that
$$H^i(V, \O_V(K_V + \rndup{M}) \otimes {\Cal J}_M) = 0 \ \text{ for }
i > 0. \tag 6.7.1$$  [Let $f : W \lra V$ be an embedded resolution of
supp$\{ M \}$. One can write $$K_W +  \rndup{f^* M} = f^*(K_V +
\rndup{M}) + P - N,$$ where $P$ and $N$ are relatively prime effective
$f$-exceptional divisors on $W$. Then set ${\Cal J}_M = f_*\O_W(-N)$.
For (6.7.1), show first that $f_*\O_P(P) = 0$. Then argue e.g. as in
the proof of Grauert-Riemenschneider vanishing in \cite{Kol1},
Corollary 11, to prove that
$R^if_*(K_W +  \rndup{f^*M})= 0$ for $i > 0$.] Esnault and Viehweg
show that ${\Cal J}_M$ is actually independent of the resolution
chosen. This is the algebro-geometric analogue --- and a special case
of --- the multiplier ideal sheaf associated to a line bundle with a
singular metric (cf. \cite{De2} and \cite{Siu2}), and (6.7.1) is a
special case of Nadel's vanishing theorem \cite{Nad}. One thinks of
${\Cal J}_M$ as measuring the singularities of the fractional part $\{
M \}$ of $M$. Its co-support is contained in the locus of points at
which supp$\{ M \}$ fails to be a normal crossing divisor.
(Compare \cite{Kol1}, Theorem 19$^\prime$.) \qed
\endxca

\xca{Exercise 6.8}  ({\bf Singularities of Plane Curves.}) In this
exercise, $C \subset \P^2$ is a reduced plane curve of degree $d$.
\sbl (i). Let $\Sigma = \text{Sing }C$ be the singular locus of $C$,
considered as a reduced subscheme of
$\P^2$. Use Kawamata-Viehweg vanishing to prove the classical theorem
that $\Sigma$ imposes independent conditions on curves of degrees $k
\ge d-2$, i.e.
$$H^1( \P^2, \I_{\Sigma}(k)) = 0 \ \text{ for } k \ge d-2. \tag *$$
[Write $\Sigma = \{ x_1 , \dots, x_m \}$, and consider the blowing-up
$f : Y \lra \P^2$ of $\P^2$ along $\Sigma$, with $E = \sum E_i$ the
exceptional divisor, and $H$ the pull-back to $Y$ of the hyperplane
divisor on $\P^2$. Let $C^\prime \subset Y$ be the proper transform of
$C$, and put  $q_i =
\text{ mult}_{x_i}(C)$, so that $f^* C = C^\prime + \sum q_i E_i$. If
$q_i \ge 3$ for all $i$, then one can argue much as in the proof of
Proposition 6.4. In case some $q_i = 2$ one has to make small
perturbations of the divisors in question, as follows. Fix first a
small rational number $0 < \eps
\ll 1$ such that $H - \eps E$ is ample. Then choose $a < 1$ such that
$a > 2/(q_i + \eps)$ for all
$i$, and consider the $\Q$-divisor  $M = (d+1)H - a C^\prime - \sum
a(q_i + \eps) E_i$. The numerical equivalence
$$M \equiv (d+1)(1-a)H + a(H - \sum \eps E_i) + a(dH - C^\prime - \sum
q_i E_i)$$ shows that $M$ is big and nef (in fact ample). Now  apply
(6.2).] Observe that this argument shows that one can replace the
ideal $\I_{\Sigma}$ in (*) by ${\Cal J} = \I_{x_1}^{q_1-1}
\cdots
\I_{x_m}^{q_m-1}$.
\sbl (ii). Suppose now that $C$ has a certain number of simple cusps
(i.e. $z^2 = w^3$ in local analytic coordinates): let $\Delta = \{
x_1, \dots , x_{\kappa} \}$ be the set of such. Prove the theorem of
Zariski \cite{Z}, \S6, that
$$H^1(\P^2, \I_{\Delta} (k))=0 \ \text{ for } k > \frac{5}{6}d - 3.$$
It follows for example that a septic curve can have no more than ten
simple cusps. [Work on a blowing-up
$f : Y \lra \P^2$ which is an embedded resolution of $C$ over the
cusps.] It is interesting to note that Zariski proves this result by
studying the irregularity of cyclic coverings of $\P^2$ branched along
$C$. One can see the approach to  Kawamata-Viehweg vanishing via
covering constructions as a vast generalization of this idea. Building
in part on work of Esnault \cite{Es}, Sakai \cite{Sak3} has applied
these and other techniques to study the singularities of plane curves.
\qed
\endxca

\bl
\xca{Exercise 6.9} ({\bf Variant of Proposition 6.4.}) Let $L$ be an
ample line bundle on a smooth projective surface $X$, and suppose
that $Z \subset X$ is a finite set. Assume given a divisor $D \in
|kL|$ for some $k > 0$, say $D = \sum d_i D_i$, plus integers $s \ge
0$ and $q > (s+2)k$, such that mult$_x(D) \ge q$ for all $x \in Z$,
and with $(s+2)d_i \le q$ whenever $D_i$  meets $Z$. Prove that then
$H^i(X, \O_X(K_X + L) \otimes \I_Z^{s+1}) = 0$ for $i > 0$. [If $(s+2)
d_i < q$, proceed as in the proof of Proposition 6.4.  In general, as
in Exercise 6.8.(i), one can introduce small perturbations, as
follows. Let $f : Y \lra X$ be the blowing-up of
$Z$, with exceptional divisor $E$, and fix $\eps > 0$ such that $f^*L
- \eps E$ is ample. For rational numbers $0 < b \ll 1$ and $0 < 1-a
\ll 1$, consider
$$
\aligned
M &= f^* \left ( L - a \left(  \frac{s+2}{q} \right ) D
\right ) - b\eps E \\
   & \equiv \left( (1-b) - \frac{a(s+2)}{q}k \right ) f^*L + b \left (
f^*L - \eps E \right).
\endaligned
$$
If $0  < b < 1 - k(s+2)/q$ and $ 1 - a < b\eps/(s+2)$, then $M$ is
ample, and every component of $E$ has coefficient $\le -(s+2)$ in $M$.]
\qed
\endxca

\sbl
\bl
%\ni {\it Reider's Theorem Revisited.}
\subhead{ Reider's Theorem Revisited.} \endsubhead
\sbl
 Our next goal is to use Kawamata-Viehweg vanishing to give another
proof of part of Reider's theorem. Specifically, we focus on the
following statement:
\proclaim{Proposition 6.10} Let $X$ be a smooth projective surface, $x
\in X$ a fixed point,  and
$L$ a nef line bundle on $X$. Assume that $L^2 \ge 5$, and that
$$L \cdot C \ge 2 \ \text{ for all curves } C \subset X \ \text{ such
that } C \owns x.$$ Then $\O_X(K_X + L)$ has a section which doesn't
vanish at $x$.
\endproclaim
\bl
\ni The plan is this: the hypothesis $L^2 \ge 5$ allows us to
construct a divisor $D \in |kL|$ ($k
\gg 0$) with high multiplicity at $x$. If $D$ has small multiplicity
at nearby points, then Proposition 6.4 applies. In the contrary case,
$D$ contains a component $D_0$ appearing with high multiplicity. The
philosophy of the KRS method is that vanishing still gives a useful
surjectivity statement. The problem is then reduced to producing a
section on $D_0$, and this is attacked using the lower bound on $L
\cdot D_0$. We remark that Sakai \cite{Sak2}, using some ideas of
Serrano,  has given a cohomological proof of Reider's theorem via
Miyaoka's vanishing theorem for Zariski decompositions of linear
series on surfaces. While the present approach has the advantage of
using mainly  general techniques from higher dimensional geometry, it
would be interesting to understand more clearly than one does at the
moment the precise connections between the two proofs. (One can see the
construction below as the first step in producing the Zariski
decomposition of the relevant linear series.)

\demo{Proof of Proposition 6.10}  The first step is to show that for
$k \gg 0$ there exists a divisor
$$ D \in |kL| \ \text{ with } \  q =_{\text{def}}  \text{ mult}_x D >
2k.  $$ This is an elementary parameter count. In fact, consider the
exact sequence
$$0 \lra H^0(\O_X(kL) \otimes \I_x^{2k + 1} ) \lra H^0( \O_X(kL)) \lra
H^0( \O_X(kL) \otimes
\O_X / \I_x^{2k + 1}). $$ We are looking for the divisor of a non-zero
section in the group on the left. Since $L$ is nef and
$L^2 \ge 5$, a standard argument using Riemann-Roch (cf. \cite{H2},
V.1.8, or \cite{SS}, p. 146) shows that:
$$ h^0(\O_X(kL)) = \frac{k^2 L^2}{2} + o(k^2) \ge \frac{5}{2} k^2 +
o(k^2). $$ On the other hand,
$$h^0(\O_X(kL) \otimes \O_X / \I_x^{2k + 1}) = \binom{2k + 2}{2} =
\frac{4k^2}{2} + o(k^2).$$ Hence $ H^0( \O_X(kL) \otimes \I_x^{2k +
1}) \ne 0$ for $k \gg 0$, as required.

\bl Fix such  a divisor $D$, and write
$$ D = \sum d_i D_i + F, $$ where $\{ D_i \}$ are the components of
$D$ passing through $x$, and $F$ is the effective (or zero) divisor
consisting of those components of $D$ disjoint from $x$. If $q > 2d_i$
for all
$i$, then Proposition 6.4 applies, and we are done.

\bl Assume next that $q < 2d_i$ for some $i$. Note that
$$ q = \sum d_i \text{ mult}_x(D_i), $$ so in the first place $q \ge
d_i$ for all $i$. It follows also  that there is a {\it unique}
component
$D_i$ --- call it $D_0$ --- of maximal mutiplicity $d_0 > q/2$.
Furthermore, $D_0$ is necessarily smooth at $x$. Consider now the
$\Q$-divisor $M = L - \frac{1}{d_0}D$ on $X$. Then
$$K_X + \rndup{M} = K_X + L - D_0 -N,$$ where $N = [\frac{1}{d_0} F]$
is an effective divisor supported away from $x$. We have the numerical
equivalence
$$ M \equiv L - \left ( \frac{1}{d_0} \right ) kL \equiv \left ( 1 -
\frac{k}{d_0} \right ) L,$$ and since $d_0 > q/2 > k$, it follows that
$M$ is nef and big. Keeping in mind that we don't need to worry about
normal crossings, Theorem 6.2 gives the vanishing $H^1(X, \O_X(K_X + L
- D_0 -N))= 0$. Therefore the restriction
$$H^0(X, \O_X(K_X + L - N)) \lra H^0(D_0, \O_{X}(K_X + L - N)
\restr{D_0} ) \tag 6.11$$ is surjective.

\bl Observe next that it is enough to show that $\O_{X}(K_X + L - N)
\restr{D_0}$ has a section $\bar t$ which does not vanish at $x$. For
then thanks to the surjectivity of (6.11), $\bar t$ lifts to a section
$t \in H^0(X, \O_X(K_X + L - N))$ with $t(x) \ne 0$. Since $x \notin
\text{ supp } N$, $t$   gives rise to a section $s \in H^0(X, \O_X(K_X
+ L))$ which is non-zero at $x$, as desired.

\bl The existence of the required section $\bar t$ will in turn follow
if we verify:
$$ (L - D_0 - N) \cdot D_0 > 1 . \tag 6.12$$ For then, since $(L - D_0
- N) \cdot D_0$ is in any event an integer, it will follow that the
restriction $\O_X(K_X + L - N) \restr{D_0}$ is of the form
$\O_{D_0}(K_{D_0} + B)$ for some line bundle $B$ of degree $\ge 2$,
and hence is free at $x$ thanks to Theorem 1.1. (Note that the
possible singularities of $D_0$ don't cause any problems here since
$D_0$ is Gorenstein and $x$ is a smooth point.)  As for (6.12), note
that $\rndup{M} - M = \sum_{i \ge 1} \frac{d_i}{d_0} D_i + \Delta$,
where
$\Delta = \rndup{\frac{1}{d_0}F} - \frac{1}{d_0}F$ is an effective
divisor which meets $D_0$  properly. Thus
$$ L - D_0 - N = \rndup{M}  \equiv M + \sum_{i \ge 1} \frac{d_i}{d_0}
D_i +\Delta.
$$  But $L \cdot D_0 \ge 2$ by assumption, $D_i \cdot D_0 \ge
mult_x(D_i)$ for $i
\ge 1$, and $M \equiv (1 - \frac{k}{d_0})L$. Thus:
$$
\aligned (L - D_0 - N ) \cdot D_0 &= \left ( 1 - \frac{k}{d_0} \right
) L \cdot D_0 + \sum_{i \ge 1}
\left ( \frac{d_i}{d_0} \right ) D_i \cdot D_0 + \Delta \cdot D_0 \\
&\ge 2 \left ( 1 - \frac{k}{d_0} \right ) + \sum_{i \ge 1}
\frac{d_i}{d_0} \text{ mult}_xD_i
 \\  &= 2 \left ( 1 -\frac{k}{d_0} \right ) + \left (
\frac{q-d_0}{d_0}  \right ) = 1 + \left (
\frac{q - 2k}{d_0} \right ).
\endaligned
$$ Recalling that $q > 2k$, (6.12) follows.

\bl It remains to treat the possibility that $q = 2 d_i$ for some $i$.
When $L$ is ample one can invoke Exercise 6.9 (i.e. introduce small
perturbations of the divisors in question and argue as in Proposition
6.4). In general,  apply Exercise 6.13  below. We leave
details to the reader.
\qed
\enddemo
\bl
\demo{Remark} This is the model of the argument used in \cite{EL2},
\cite{ELM} and \cite{Fuj2} to study global generation of  adjoint
linear series on threefolds. Given an ample line bundle $L$ on a smooth
projective threefold $V$, plus a point $x \in V$,  one starts by taking
a divisor $D \in |kL|$ $(k \gg 0)$ with high multiplicity at $x$.
If  $D$ has an (almost) isolated singularity at $x$, then one
concludes at once as in Remark 6.5.  Otherwise --- very roughly
speaking ---  the KRS approach reduces one to producing a section on
the ``most singular locus" $Z$ of $D$, which in the case at
hand is either a curve or a component appearing with particularly high
multiplicity in $D$. Then one applies Reider-type statements for
$\Q$-divisors on $Z$.  However there are considerable technical
difficulties stemming in part from the fact that one has to start by
passing to an embedded resolution of
$D$, since already in dimension three one can't ignore the normal
crossing hypothesis in (6.2).  Unfortunately, these problems have
so far blocked the possibility of extending  the argument to dimensions
four or more.
\enddemo
\bl

\xca{Exercise 6.13}  Let $R$ be a big and nef $\Q$-divisor on a smooth
projective surface $X$, and let $E_1, \dots, E_k$ be distinct
irreducible curves on $X$ which do not appear in the fractional part
of $R$. Assume that $R \cdot E_i > 0$ for all $1 \le i \le k$. Prove
that  then
$$H^i(X, \O_X(K_X + \rndup{ \ R \ } + E_1 + \dots + E_k )) = 0 $$
for $i > 0$. [Use induction on $k$.]
\qed \endxca

\bl
\xca{Exercise 6.14} ({\bf Reider-type Theorems for $\Q$-divisors},
\cite {EL2},  \S2.) The proof of Proposition 6.9 can be extended to
give analogous statements for (round-ups of) $\Q$-divisors. Let $X$ be
a smooth projective surface.
\sbl (i).  Suppose that $M$ is an ample $\Q$-divisor on $X$ such
that
$$M^2 > 4 \ \ \text{and} \  \ M \cdot C \ge 2\text{ for all
irreducible curves $C \subset X$}.$$ Show that then $\O_X(K_X +
\rndup{M})$ is globally generated. [Fix $x \in X$ and then a divisor
$D \in |kM|$ for sufficiently divisible $k \gg 0$ such that $q =
\text{mult}_x(D)  > 2k$. Write $D = \sum d_i D_i + F$ and $M =
\rndup{M} - \sum a_i D_i - G$ ($0 \le a_i <1$) where the $D_i$ all
pass through $x$, and $F, G$ are effective divisors disjoint from $x$.
Let $p = \text{mult}_x(\rndup{M} - M) = \sum a_i \ \text{mult}_x D_i$.
If $(2 - p)/q < (1 -a_i)/d_i$ for all $i$, blow up $x$ and argue as in
the proof of Proposition 6.4. If the reverse inequality holds for some
$i$, adapt the proof given for Proposition 6.10.]
\sbl (ii). For applications, it is useful to have a relative
statement.  Thus consider a surjective map
$h : X \lra X_0$ where $X_0$ is a complete irreducible surface. Let
$M$ be a big and nef
$\Q$-divisor on $X$, and fix a point $x_0 \in X_0$.  Suppose that
$\beta_1, \beta_2 > 0$ are positive rational numbers such that
$$
\align M^2 &> (\beta_2)^2 \\ M \cdot C &\ge \ \beta_1 \ \ \forall
\text{ curves } C \subset X \text{ s.t. $h(C)$ is a curve through }
x_0.
\endalign
$$ Assume that $ \O_X(K_X + \rndup{M}) | \Gamma \cong \O_{\Gamma}$ for
every effective divisor $\Gamma
\subset X$ such that $h(\Gamma) = x_0$. Suppose also that
$$ \beta_2 \ge 2, \ \ \ \beta_1 \left ( 1 - \frac{2}{\beta_2} \right )
\ge 1.$$ Show that then $ \O_X(K_X + \rndup{M}) $ has a section which
is non-vanishing at some point $x \in h^{-1}(x_0)$. [See
\cite{EL2}, Theorem 2.3.]  \qed
\endxca

\bl
\xca{Exercise 6.15} ({\bf Reider-type Theorem for Normal Surfaces},
\cite{ELM}, \S1, compare \cite{Sak2}.) Let $S$ be a complete normal
surface. Recall that Mumford  has defined a $\Q$-valued
intersection product for Weil divisors on $S$. [In brief, let $f : T
\lra S$ be a resolution of S. Given a Weil divisor $D$ on $S$, let
$D^\prime$ denote its proper transform  on $T$. Then there is a unique
$f$-exceptional $\Q$-divisor $\Delta$ on $T$ such that $$(D^\prime +
\Delta) \cdot E = 0$$ for all $f$-exceptional divisors $E$ on $T$.
Mumford first defines $f^* D = D^\prime + \Delta$. Given two
Weil-divisors $D_1 , D_2$ on $S$, one then sets $D_1 \cdot D_2 = f^*
D_1 \cdot f^* D_2 \in \Q$.]  In particular the usual definition for
nefness makes sense for Weil divisors on $S$. Similar considerations
of course hold for $\Q$-divisors on $S$. Recall also that the
canonical divisor $K_S$ exists as a Weil divisor (class) on $S$.
\sbl
	Now suppose that $M$ is a nef $\Q$-divisor on $S$, and let $\beta_1,
\beta_2 > 0$ be rational numbers such that
$$M ^2 > (\beta_2)^2 \ \ , \  \ M \cdot C \ge \beta_1 \ \  \forall  \
\text{curves } C \subset S.$$ Assume that $K_S + \rndup{M}$ is
Cartier, and that
$$\beta_2 \ge 2 \ \ , \ \ \beta_1 \left( 1 - \frac{2}{\beta_2} \right
) \ge 1. \tag * $$ Prove that then $\O_S(K_S + \rndup{M})$ is globally
generated. The inequalities (*) are satisfied for example if $M^2 >
16$ and $M \cdot C \ge 2$ for all curves $C$. We don't know whether
the conclusion holds assuming only that $M^2 > 4$ and $M \cdot C \ge
2$. [Let $f : T \lra S$ be the minimal resolution of $S$, and define
a divisor $W$ on $T$ by the equation
$K_T + \rndup{f^*M} = f^*(K_S + \rndup{M}) - W$. Prove that $W$ is
integral and exceptional. Show that every component of $f^*\rndup{M} -
\rndup{f^*M}$ has coefficient $> -1$, and deduce that $W$ is
effective (or zero). Now apply Exercise 6.14.(ii).]
\qed
\endxca

\bl
\xca{Exercise 6.16} ({\bf Very ampleness for $\Q$-divisors},
\cite{Mas}.) Let $X$ be a smooth projective surface, and let $M$ be
a $\Q$-divisor on $X$. Assume that
$$ M^2 > 18  \ \ \text{ and } \ \ M \cdot C \ge 3 \
\text{for all effective curves } C \subset X.$$
Prove that then $\O_X(K_X + \rndup{M})$ is very ample. [The linear
series in question is globally generated thanks to Exercise
6.14.(i), and the first step is to argue that given distinct points
$x , y \in X$, one can find a divisor $\Lambda \in |K_X +
\rndup{M}|$ passing through one but not the other of the points. To
this end, choose $D\in |kM|$, for $k$ very large and divisible,
such that  $q_x = \text{mult}_x(D) > 3k$ and $q_y = \text{mult}_y(D) > 3k$.
Write $\Delta =_{\text{def}} \rndup{M} - M = \sum a_j D_j$ and $D =   \sum r_j
D_j $, set $\mu_x = \text{mult}_x(\Delta)$, $\mu_y = \text{mult}_y(\Delta)$,
and
put
$$c_x = \min \left \{ \frac{2 - \mu_x}{q_x} , \frac{1 - a_j}{r_j} \biggm | D_j
\owns x \right \},$$
with $c_y$ defined similarly. Define $c = \max \{ c_x ,
c_y \}$, and consider the $\Q$-divisor $M - cD = \rndup{M} - \Delta - cD$.
Assuming for concreteness that $c = c_x$, one argues in cases according to
whether the value of $c_x$ is determined by the multiplicities at $x$, or
whether it is accounted for by a  component $D_j$ -- say $D_0$ -- of
$D$ passing through $x$. In the latter instance, one further distinguishes
between whether or not $D_0$ goes also through $y$. The argument for separating
tangent directions at $x$ is similar, but requires a little more thought. See
\cite{Mas} for details, as well as for a statement involving  more
general numerical hypotheses in the spirit of Exercise 6.14.(ii).]
\qed
\endxca

\bl
\xca{Exercise 6.17} ({\bf Shokurov's Non-Vanishing Theorem},
\cite{Sho}, cf. \cite{CKM}, Lecture 13.) In this exercise, we outline
the proof of a ``toy" version of Shokurov's Non-Vanishing Theorem. In
the case of surfaces treated here, the statement is a  consequence of
Riemann-Roch. However the argument we indicate, suitably modified,
also works in higher dimensions. (As the reader will note, it was the
inspiration for the proof of Proposition 6.10 as well as
\cite{EL2}.) The result we aim for is the following:
\proclaim{Theorem} Let $X$ be a smooth projective
surface, and let $L$ be a nef line bundle on $X$. Assume that $$A
=_{\text{def}}L - K_X \ \ \text{is ample.}$$ Then $H^0(X, \O_X(mL))
\ne 0$ for all $m \gg 0$. \endproclaim

\sbl (i).  If $L$ is numerically trivial, show that the statement
follows directly from a computation of
$\chi(X,
\O_X(mL))$. Hence we may  assume that $L$ is not numerically trivial.
Deduce  that then the intersection number
$(pL -K_X)^2$ is an increasing function of $p$.
\sbl  (ii). Fix a general point $x \in X$. Show that one can take $p$
sufficiently large so that for
$k
\gg 0$ there exists a divisor $D \in |k(pL - K_X)|$ such that $q
=_{\text{def}}  \text{mult}_x(D) > 2k$. Write $D = \sum d_i D_i$, and
set $d = \max \{ d_i \}$. Given a rational number $c > 0$, consider the
$\Q$-divisor $N(m, c) = mL - K_X - cD$. The numerical equivalence
$$ mL - K_X - cD \equiv \left ( 1 - ck \right )
\left (  (L - K_X) + (p-1)L
\right ) + (m - p) L
$$
shows that $N(m,c)$ is ample if $m \ge p$ and $ck < 1$.

\sbl (iii). If $q > 2d$, let $f : Y \lra X$ be the blowing-up of $X$ at
$x$. Argue as in Proposition 6.4 with the $\Q$-divisor $f^* N(m,
\frac{2}{q})$ on $Y$ to deduce the desired non-vanishing.
\sbl
(iv). Assume $q \le 2d$. Let $D_0$
be a component of $D$ of maximal multiplicity
$d_0 = d$, and consider the $\Q$-divisor $N(m, \frac{1}{d_0})$. Deduce
from Vanishing and Exercise 6.13 that if $m \ge p$, then the
restriction map
$$H^0(X, \O_X(mL)) \lra H^0(D_0, \O_{D_0}(mL)) \tag *$$ is surjective.

\sbl (v). In the situation of (iv), show that if $m \ge p$, then
$$\O_{D_0}(mL) = \O_{D_0}(K_{D_0} + A_0) \tag ** $$ for some ample
line bundle $A_0$ on $D_0$. Conclude that $H^0(D_0, \O_{D_0}(mL)) \ne
0$ for $m \ge p$,  and show that the theorem follows. [If $p_a(D_0)
\ge 1$ the required non-vanishing follows from Riemann-Roch. When
$D_0$ is rational, note that the bundle on the right in (**) can't
have degree $-1$ since its degree is divisible by $m \gg 0$.]
\qed \endxca

\bl
\xca{Exercise 6.18} ({\bf Kawamata's Basepoint-free Theorem},
\cite{K2}, \cite{K3}, cf. \cite{CKM}  Lecture 10.) This exercise is
devoted to a stripped-down version on surfaces of a fundamental
theorem of Kawamata:
\proclaim{Theorem} Let $X$ be a smooth
surface, and let $L$ be a nef  line bundle on $X$ such that $A
=_{\text{def}}L - K_X $ is ample. Then the linear series $|mL|$ is
free for all $m \gg 0$. \endproclaim

\sbl
 (i). Observe to begin with it follows from the
Non-Vanishing Theorem (Exercise 6.17) that $H^0(X, \O_X(mL)) \ne 0$
for all $m \gg 0$. Let $B(m)$ denote the reduced base locus of $|mL|$.
Noting that $B(c^a) \subseteq B(c^b)$ whenever $a >b$ and $c \gg 0$,
show that the sequence of subsets $B(c^n)$ stabilizes for $n \gg 0$.
Denoting the limit by $B(c)$, show that it is enough to prove that
$B(c) = \emptyset$ for every $c \gg 0$.

\sbl  (ii). Suppose first that for some $m \ge 1$ the linear series
$|mL|$ has only isolated base-points. Show that then $|(3m+a)L|$ is
free for all $a \ge 1$. [If the base-locus of $|mL|$ consists of a
finite set $Z \subset X$, then there exists a reduced divisor $D \in
|3mL|$ with mult$_x(D) \ge 3$ for every
$x \in Z$. Let $f : Y  \lra X$ be the blowing-up of $Z$, and consider
the pull-back to $Y$ of the $\Q$-divisor
$$ \left (3m + a \right ) L - \frac{2}{3}D \equiv K_X + A + (m+ a -
1)L. $$
Use vanishing to deduce that $H^1(X, \O_X((3m + a)L \otimes \I_Z)) =
0$. ]
\sbl Next, fix some large integer $p$ and write
$$|pL| = |M| + \sum r_i F_i ,$$ where $F = \sum r_i F_i$ is the fixed
divisor of $|pL|$, and $|M|$ is the moving part of $|pL|$, so that
$|M|$ has at most isolated fixed points. Assuming $F \ne 0$, we will
argue that if $F_0$ is a component appearing with maximal multiplicity
in $F$, then $F_0$ is no longer a fixed component of $|mL|$ for $m \gg
0$. In view of (i) and (ii), the theorem will follow. Turning to
details,  let $D
\in |pL|$ be a general divisor. Then
$D = D_1 + \sum  r_i F_i$ where $D_1$ is reduced. After re-indexing if
necessary we can suppose that
$r_0 = \max{r_i}$.  Consider the $\Q$-divisor
$N(m) = mL -  K_X - \frac{1}{r_0}D$. Since $L$ is nef, the numerical
equivalence
$$N(m) \equiv \left ( m - 1 - \frac{p}{r_0} \right ) L + A$$ shows
that  $N(m)$ is ample if $m > p+1$.
\sbl (iii). Keeping the notation just introduced, use Vanishing and
Exercise 6.13  to show that the restriction map
$$H^0(X, \O_X(mL)) \lra H^0(F_0, \O_{F_0}(mL)) \tag *$$  is surjective
when $m > p+1$. Arguing as in Exercise 6.17.(v), show that the group
on the right in (*) is non-vanishing.  Conclude that $F_0$ is not in
the base locus of $|mL|$ for $m > p + 1$. \qed
\endxca

\bl
\xca{Exercise 6.19} ({\bf Pluricanonical Series}.)   With some extra
work, one can show that the conclusion of Kawamata's Basepoint-free
theorem (Exercise 6.18) remains valid assuming only that $L - K_X$ is
nef and big (cf. \cite{CKM}, Lecture 10). We grant this fact here.
\sbl (i). Show that if
$X$ is a minimal surface of general type, then $\O_X(mK_X)$ is
globally generated for all $m
\gg 0$. Deduce that there exists a normal surface $X_0$, a surjective
birational morphism $h : X
\lra X_0$, plus an ample line bundle $L_0$ on $X_0$ such that
$\O_X(K_X) = h^* L_0$. ($X_0$ is called the {\it canonical model} of
$X$. Compare  \cite{EV2}, \S 7.1, and \cite{W}.)
\sbl
(ii). With $X$ as in (i), use Exercise 6.14.(ii) to prove  that
$\O_X(mK_X)$ is globally generated for $m \ge 5$.
\sbl
\ni This is the model of the arguments used to study pluricanonical
series on threefolds in \cite{EL2} and \cite{ELM}. Note however that
the bound in (ii) is slightly less than optimal.
\qed
\endxca

\sbl
\bl
%\ni {\it Bogomolov's Theorem.}
\subhead {Bogomolov's Theorem.} \endsubhead
\sbl Fern\'andez del Busto \cite{FdB1} has shown that the sort of
argument used to prove Propositions 6.4 and 6.10 leads to a new
approach to Bogomolov's Instability Theorem 4.2. 	This completes in a
very nice way the circle of ideas linking linear series on surfaces,
vector bundles and vanishing theorems, and we present here an outline
of his proof.  We shall content ourselves with explaining the
principal steps. Therefore  we'll make some simplifying assumptions,
and relegate some of the numerical calculations to  exercises.

\bl
To begin with, recall the statement:
\proclaim{Bogomolov's Instability Theorem} Let $E$ be a rank two
vector bundle on a smooth projective surface $X$. Assume that
$$ c_1(E)^2 - 4c_2(E) > 0.  $$ Then there is a saturated invertible
subsheaf $A \hookrightarrow E$ such that if $L = det \ E$, then
$$(2A - L)^2 > 0 \ \ \text{and} \ \ (2A - L) \cdot H > 0 \ \text{for
some ample divisor } H.$$
\endproclaim
\ni By a saturated subsheaf we mean here that the vector bundle map $A
\lra E$ vanishes only on a finite set. (See Definition 4.1.)
\bl We start with some preliminary reductions and remarks. Note first
that since Bogomolov's theorem is invariant under twisting, we are
free to  tensor $E$ by a high multiple of an ample line bundle.
Therefore we may  assume that $E$ is globally generated, that its
determinant is ample, and that $c_2(E) > 0$. Let $s \in \Gamma(X,E)$
be a general section. Then  the zero-scheme $Z = Z(s)$ is reduced, and
the Koszul complex (3.6) determined by $s$ realizes $E$ as an extension
$$0  \lra \O_X \overset{ \cdot s} \to \lra E \lra L \otimes \I_{ Z}
\lra 0 . \tag 6.20 $$ Here $L = \text{det}(E)$ is an ample line bundle
on $X$, and $c_2(E) = \# Z > 0$.

\bl The strategy of the proof is now very simple. The numerical
hypothesis $c_1^2 > 4c_2$ guarantees the existence of a divisor
$D \in |kL|$ ($k \gg 0$) with high multiplicity at the points of $Z$.
If $D$ has small multiplicity away from $Z$, then the argument of
Proposition 6.4 would yield the vanishing of $H^1(X,
\O_X(K_X + L) \otimes \I_Z)$. But since $E$ is locally free, Theorem
3.13 shows that the group in question is non-zero. Therefore, as in
the proof of (6.10), $D$ must contain some distinguished components
appearing with high multiplicity. One uses these components to
construct a divisor $\Gamma
\subset X$ passing through $Z$. Vanishing will imply that the
inclusion $\O_X(L - \Gamma)
\hookrightarrow L \otimes
\I_Z$ lifts to an embedding  $\O_X(L - \Gamma) \subset E$, and then
 one argues that the saturation of this subsheaf destabilizes
$E$. One new feature of this argument is that while in Proposition
6.10 we worked with an arbitrary divisor $D$, here we must be careful
to choose
$s$ and $D$ rather generally.
\bl
	Turning to the details, the first point is:
\proclaim{Lemma 6.21} For any section $s \in \Gamma(X, E)$ as above,
with reduced zero-scheme
$Z  = Z(s)$, there exists a divisor
$D \in | kL | $ $(k \gg 0)$ such that mult$_x(D) > 2k$ for every $ x
\in Z(s)$. Moreover, by choosing
$s$ and $D$ sufficiently generally, we can assume that $D$ satisfies
the following
\sbl
\ni {\smc Uniform Multiplicity Property  (UMP). }{ For any rational
number $\delta > 0$, the multiplicity of the integer part
$[ \delta D ]$ is the same at every point of $Z(s)$.}
\endproclaim
\demo{Proof} For the first assertion, one counts dimensions much as in
the proof of Proposition 6.10. Specifically, by Riemann-Roch:
$$h^0(X, \O_X(kL)) =  \frac{k^2}{2} L^2  + o(k^2) =
\frac{k^2}{2}c_1(E)^2 + o(k^2).$$ On the other hand,  the number of
conditions required to impose multiplicity
$\ge 2k + 1$ at each point of $Z$ is:
$$ \binom{2k+2}{2} (\# Z) =  \binom{2k+2}{2} c_2(E) =
\frac{4k^2}{2}c_2(E) + o(k^2).$$ Since
$c_1(E)^2 > 4c_2(E)$ by hypothesis, for $k\gg 0$ the required divisor
will exist.

\bl
\def \ZZ{\Cal {Z}}
\def \DD{ \Cal {D}} Turning to the second assertion, let $S \subset \P
H^0(E)$  denote the open subset   parametrizing sections of $E$ with
finite zero-schemes, and consider  the incidence correspondence
$$X \times S \supset \ZZ = \left\{ (x, [s]) \bigm | x \in Z(s) \right
\}.$$ The assumption that
$E$ is globally generated implies that $\ZZ$ is an open subset of a
projective bundle over $X$, and hence is irreducible. Note also that
it is finite over $S$, and in particular dim$ (\ZZ)  =
\text{dim}(S)$.  By an evident globalization of the parameter count
just made, one can construct a divisor
$X \times S \supset \DD \supset \ZZ$ such that the fibre $D_s \subset
X$ of $\DD $ over $s \in S$ is a divisor in the linear series $|kL|$
having multiplicity $> 2k$ at every point of $Z_s = Z(s)$.
\bl
 To establish the UMP, it is enough to show that  given an integer $p
> 0$, if $[\delta D_s]$ has multiplicity $\ge p$ at {\it some} point
$x \in Z(s)$ for general $s \in S$, then it has multiplicity $\ge p$
at {\it every} point $x \in Z(s)$. Consider to this end the
Zariski-closed set
$$\ZZ \supset \ZZ_p = \left \{ (x,[s]) \in \ZZ \bigm |
\text{mult}_x([\delta \DD]_s) \ge p \right
\}.$$ By assumption $\ZZ_p$ dominates $S$, and hence dim$(\ZZ_p) \ge
\text{dim}(S) = \text{dim}(\ZZ)$.
$\ZZ$ being irreducible, it follows that $\ZZ_p = \ZZ$. Noting that
$\text{mult}_x([\delta \DD]_s) =
\text{mult}_x([\delta D_s])$ for general $s$, this means exactly that
$[\delta D_s]$ has multiplicity $\ge p$ for every $x \in Z(s)$. The
Lemma follows. \qed
\enddemo
\bl
\xca{Remark} Observe that the second assertion of the Lemma is
essentially a monodromy argument in the spirit of \cite{Har}.
Specifically, fix a general ``reference section" $s_0 \in
\Gamma(X, E)$ with finite reduced zero-scheme $Z_0$. Letting $s$ vary
over the open subset $U
\subset S$ parametrizing sections with reduced zero-loci, one obtains
a monodromy action on the points of $Z_0$. The irreducibility of $\ZZ$
implies that this action is transitive. Since the divisors $D_s$ also
vary with $s$, it follows that for general $s \in U$, $D_s$ cannot be
used to distinguish among the points of $Z(s)$. Therefore $D_s$ must
have the same multiplicity at every point of $Z(s)$, and similarly for
$[\delta D_s]$.
\endxca
\bl

Returning to the proof of Bogomolov's theorem, fix $s$ and $D \in
|kL|$ ($k \gg 0$), with $D \supset Z = Z(s)$, as in Lemma 6.21. In
particular, we suppose that $D$ satisfies the Uniform Multiplicity
Property. Therefore
$D$ has the same multiplicity at every point of $Z$, say
$$ q = \text{mult}_x(D) \ \text{ for every } x \in Z.$$ Write
$ D = \sum d_i D_i$. We henceforth {\it make the simplifying
assumption} that every component
$D_i$ of $D$ meets $Z$. (See Exercise 6.24  for the general case.)
Then set $$d = \max \{ d_i \}.$$

\bl We assert next that $2d > q$. In fact, suppose to the contrary
that $2d \le q$. Then an evident generalization of Proposition 6.4
(i.e. Exercise 6.9) implies that
$$H^1(X,\O_X(K_X + L) \otimes \I_Z) = 0. \tag * $$
 On the other hand, since $E$ is locally free and $Z \ne \emptyset$,
the extension class of the Koszul complex (6.20) must be non-zero.
Therefore
$$\text{Ext}^1(L \otimes \I_Z, \O_X) = H^1(X, \O_X(K_X + L) \otimes
\I_Z)^* \ne 0,$$ a contradiction. (Compare (3.12), (3.13).) Hence $2d
> q$, as claimed.

\bl Consider now the divisor
$$D_0  = \left [ \sum \frac{d_i}{d} D_i \right ] . $$ Then $D_0$ is
reduced, and $D_0$ meets $Z$ thanks to our assumption that every
component $D_i$ passes through at least one point of $Z$. It follows
from the (UMP) that in fact $Z \subset D_0$. Furthermore, $D_0$ is
smooth at every point of $Z$. $D_0$ plays the role of the curve
$\Gamma$ appearing in the overview of the proof given above.

\bl
	Since
$$\frac{1}{d} < \frac{2}{q} < \frac{1}{k},$$  the $\Q$-divisor $L -
\frac{1}{d}D$ is ample. Kawamata-Viehweg vanishing (6.2), along with
Sakai's Lemma (Ex. 6.6), therefore applies to show  that
$$H^1(X, \O_X(K_X + L - D_0) ) = \text{Ext}^1( \O_X(L - D_0), \O_X)^*
= 0. \tag **
 $$
Now we can pull back the given extension (6.20) under the natural
inclusion $\O_X(L - D_0) \subset \O_X(L) \otimes \I_Z$ to obtain an
extension of $\O_X(L - D_0)$ by $\O_X$. But the vanishing (**)
implies that the latter extension splits. Therefore the
injection $\O_X(L - D_0) \subset \O_X(L) \otimes \I_Z$ lifts to a sheaf
monomorphism $\O_X(L - D_0) \hookrightarrow E$:
$$
\CD @. @. \O_X(L - D_0) @= \O_X(L - D_0) \\ @. @.  @VVV @VV{\cdot
D_0}V \\ 0 @>>> \O_X @>{\cdot s}>> E @>>> \O_X(L) \otimes \I_Z @>>> 0.
\endCD
$$  In other words, we have produced an ``unexpected" invertible
subsheaf of $E$. To complete the proof, Fern\'andez del Busto shows
that $\O_X(L - D_0)$ is a saturated subsheaf of $E$, and that it
destabilizes $E$. The destabilization is equivalent to proving the
inequalities:
$$
\aligned (L - 2D_0)^2 > 0 \\ (L - 2D_0) \cdot L > 0,
\endaligned \tag 6.22
$$
 and we outline the required calculations in the following Exercise.
\bl
\xca{Exercise 6.23} Keeping the notations and assumptions made above,
we sketch the verification of (6.22).
\sbl (i).  Suppose given a decomposition $Z = Z_1 \coprod Z_2$ of $Z$
into two disjoint subsets, with
$Z_1 \ne \emptyset$. Prove that $Z_1$ cannot impose independent
conditions on the sections of
$\O_X(K_X + L)$ vanishing on $Z_2$, i.e. that the evaluation
homomorphism
$$ H^0(X, \O_X(K_X +  L) \otimes \I_{Z_2}) \lra H^0(Z_1, \O_{Z_1}(K_X
+ L ))$$ cannot be surjective. [Use Theorem 3.13.]
\sbl (ii). Prove that $$(L - D_0) \cdot D_0 \le \# Z.$$ [Each point of
$Z$ lies on exactly one component of $D_0$, and every component of
$D_0$ contains at least one point of $Z$. So it suffices to show that
if $D^\prime$ is a component of $D_0$, and
$Z^\prime = D^\prime \cap Z$, then $(L - D_0) \cdot D^\prime \le \#
Z^\prime$. If on the contrary $(L - D_0) \cdot D^\prime > \text{deg }
Z^\prime$, then the points of $Z^\prime$ impose independent conditions
on the linear series $|\O_{D^\prime}(K_{D^\prime} + L - D_0)|$. Use
(**) and
 (i) to arrive at a contradiction.]
\sbl (iii). Show that
$$(L - D_0) \cdot D_0 \ge \left ( 1 - \frac{k}{d} \right ) L \cdot D_0
+ (\# Z) \left (
\frac{q}{d} - 1 \right ).$$ [$L - D_0 \equiv (1 - \frac{k}{d})L +
D^*$, where $D^* = \sum \frac{d_i}{d}D_i - D_0$ is an effective
$\Q$-divisor which does not contain any components in common
with $D_0$. Estimate $D^* \cdot D_0$ as in the proof of (6.12).]
\sbl (iv). Show that $$\# Z > \frac{1}{2} L \cdot D_0.$$ [Combine (ii)
and (iii).]
\sbl (v). Prove that $\O_X(L - D_0)$ is a saturated subsheaf of $E$,
and establish the inequalities (6.22). [For the saturation, use the
fact that  every component of
$D_0$ contains a point of $Z$ lying only on that component .  The
first inequality in (6.22) then follows formally from the hypothesis
$c_1(E)^2 > 4c_2(E)$ by a suitable computation of $c_2(E)$, much as in
(4.8). As for the second, one has $L^2 > 4 (\# Z) > 2 L \cdot
D_0$.]
\qed
\endxca

\bl
\xca{Exercise 6.24}  Eliminate from the argument above the simplifying
assumption that every component of $D$ meets $Z$. [Take $D$ satisfying
Lemma 6.21, and write $D = \sum d_i D_i + F$, where every
$D_i$ meets $Z$ and $F$ is disjoint from $Z$. Defining $q$ and $d$  as
above, one still has $2d > q$. Choose a minimal subdivisor $\Delta
\subseteq [ \frac{1}{d} F]$ for which one has the vanishing
$$H^1(X, \O_X(K_X + L - D_0 -\Delta)) = 0,$$ where  $D_0 = \sum
\frac{d_i}{d} D_i$. Set $\Gamma = D_0 + \Delta$. Then replace $D_0$ in
the argument above with $\Gamma$. See \cite{FdB1} for details.]
\qed
\endxca

\bl
\bl

\head {\S 7. Algebro-Geometric Analogue of Demailly's Approach }
\endhead
\bl

The method of Kawamata-Reid-Shokurov discussed in the previous
section is essentially inductive in nature. To produce a section
of $\O_X(K_X + L)$, for example, one starts by constructing a
divisor $D \in |kL|$ for $k \gg 0$ with large multiplicity at a
given point $x \in X$. If $D$ has small multiplicity at
neighboring points $y \in X$ --- which is the ``good" case ---
then one gets a vanishing which directly yields the required
section (Proposition 6.4). If on the contrary $D$ has a ``bad"
component $D_0$ appearing with high multiplicity, then the KRS
machine reduces the problem to constructing a section on $D_0$.
Unfortuately, already in the three dimensional situation of
\cite{EL2} and \cite{ELM}, this inductive step involves
considerable technical difficulties.

\bl Demailly's strategy in \cite{De1}  is quite different. In
effect he puts his efforts into showing that under suitable
hypotheses  one can produce a divisor $D \in |kL|$ with an
``almost isolated" singular point, thereby avoiding  any
inductions. (We're being somewhat metaphorical here:  strictly
speaking, Demailly does not work directly with divisors.)  The
argument of \cite{De1} becomes essentially analytic at one point,
but it has recently emerged in work with Ein and Nakamaye
\cite{ELN} that there are quite simple algebro-geometric
analogues of Demailly's underlying geometric ideas. In the
present section we discuss this approach in the case of surfaces,
where (as usual) the picture is particularly clear.
\bl
\sbl

  Let $X$ be a smooth projective surface, and consider a fixed
point $x \in X$ . We start with a definition:
\definition{Definition 7.1} Given a line bundle $B$ on $X$, and a
divisor  $D \in |kB|$ for some
$k > 0$, we say that $D$ has an {\it almost isolated singularity
of index} $> r$  at $x$ if
$$\text{mult}_x(D) > kr,$$  and if there exists a neighborhood $U
\owns x$ of $x$ in $X$ such that
$$\text{mult}_y D < k \ \text{for } y \in U - \{ x \} . $$
\enddefinition
\sbl \ni  In other words, we require that the $\Q$-divisor
$\frac{1}{k} D \equiv B$ have multiplicity $> r$ at $x$, and
$< 1$ in a punctured neighborhood of $x$. Similarly, one can
discuss an almost isolated singularity of index $\ge r$.
\bl
 We may then rephrase Proposition 6.4 as:
\proclaim{Proposition 7.2}  Let $X$ be a smooth surface, $L$ a
big and nef line  bundle on $X$, and $s \ge 0$ an integer. If for
some $k > 0$ there exists a divisor $D \in |kL|$ with an almost
isolated singularity of index
$> s+2$ at a point $x \in X$, then $|\O_X(K_X + L)|$ generates
$s$-jets  at $x$. \qed \endproclaim
\bl
\ni So we need to find a way to produce divisors with almost
isolated singularities.

\bl The potential obstruction to constructing such divisors is
easily described, especially on surfaces. Specifically, in the
situation of the Proposition consider the linear series
$$|\O_X(kL) \otimes  \I^{k(s+2) + 1}_x|$$   of all divisors
$D \in |kL| $ having multiplicity $ > k(s+2)$ at $x$.  Let $ F_k$
denote the fixed divisor of this linear series, and let
$\Sigma_k$ denote its moving part, so that $\Sigma_k$ has at most
isolated base points. Thus
$$| \O_X(kL) \otimes  \I^{k(s+2) + 1}_x|  = \Sigma_k + F_k,
\tag * $$  i.e. every divisor $D \in | \O_X(kL) \otimes
\I^{k(s+2) + 1}_x|$ is of the form $D = D^\prime + F_k$ for some
$D^\prime \in \Sigma_k$. Now by Bertini's theorem, a general
divisor $D^\prime \in \Sigma_k$ is reduced. Therefore a general
element $D \in | \O_X(kL) \otimes  \I^{k(s+2) + 1}_x|$ will fail
to have an almost isolated singular point at $x$ if and only if
the fixed divisor $F_k$ has components of multiplicty $\ge k$ that
pass through $x$. So the problem is prove an upper bound on the
coefficients of the components of $F_k$ for $k \gg 0$. The
strategy roughly speaking will be to use  a   lower bound  on the
dimension  of the  linear series on the left in (*) to deduce an
upper bound on the degree of $F_k$.

\bl It is useful to generalize slightly. Let $Y$ be a smooth
projective surface, and let $B$ be a big line bundle on $Y$. (In
the application, $Y$ will be a blow-up of
$X$.) Denote by $F_k$ the fixed divisor of the complete linear
series $|kB|$ and put $M_k = kB - F_k$, so that
$$| kB| = |M_k| + F_k$$ for all $k$.  The essential point is:
\proclaim{Proposition 7.3}  Assume that there exists a rational
number $\rho > 0$ such that
$$h^0(Y, \O_Y(kB)) \ge \rho \frac{k^2}{2} + o(k^2) \ \
\text{for } k \gg 0.$$ Then for $k \gg 0$:
$$M_k^2 \ge \rho k^2 + o(k^2).$$ \endproclaim

\bl
\ni To understand the intuition here, suppose it were to happen
that $M_k = k M$ for some fixed nef and big divisor
$M$. Then
$$\rho \frac{k^2}{2} + o(k^2) \le h^0(\O_Y(kB)) = h^0(\O_Y(kM)) =
\frac{k^2}{2} M^2 + o(k^2) $$ and the desired inequality follows.
The content of the Proposition is that the same inequality holds
in general.

\bl The following argument is due to Fern\'andez del Busto
\cite{FdB2}; the original proof of (7.3) was more cumbersome.
\demo{Proof of Proposition 7.3} Fix a very ample divisor $H$ on
$Y$ with the property that $K_Y + H$ is very ample, and fix a
general divisor $D \in |K_Y + H|$. Then multiplication by $D$
defines for every $k$ an inclusion of sheaves
$\O_Y(M_k) \subset \O_Y(K_Y + H + M_k)$. In particular,
$$ h^0(Y, \O_Y(K_Y + H + M_k)) \ge h^0(Y, \O_Y(M_k))  = h^0(Y,
\O_Y(kB)) \tag * $$ for all $k$. Now $M_k$ is nef, hence $H + M_k$
is ample, so by Vanishing and Riemann-Roch:
$$
\align h^0( Y, O_Y(K_Y + H + M_k)) &= \chi( Y, \O_Y(K_Y + H +
M_k)) \\ &= \frac{  \left ( K_Y + H + M_k
\right ) \cdot  \left ( H + M_k \right )   }{2}  + \chi(Y,
\O_Y) \\ &= \frac{M_k \cdot M_k}{2} + \frac{M_k \cdot \left (
K_Y  + 2H \right ) }{2} + c(Y,H),
\endalign
$$ where $c(Y,H)$ is a constant not depending on $k$. On the
other hand, since $K_Y + H$ and $H$ are very ample, there exists
a divisor $D^\prime \in |K_Y + 2H |$ which meets
$F_k$ properly. Therefore
$$ 0 \le M_k \cdot \left ( K_Y + 2H \right ) \le kB \cdot
\left ( K_Y + 2H \right ).$$ Hence
$$
 h^0( Y, O_Y(K_Y + H + M_k)) = \frac{M_k \cdot M_k}{2} +
o(k^2),$$ and the Proposition follows from (*) and the hypothesis
on $h^0(\O_Y(kB))$. \qed \enddemo
\bl
\demo{Remark} One may view Proposition 7.3 as a numerical
counterpart to the existence of a Zariski decomposition for
$B$.  It would be very interesting to know to what extent an
analogue of  (7.3) remains true on higher dimensional
varieties. The difficulty of course is that the decomposition $kB
= F_k + M_k$ of $|kB|$ into fixed and moving parts in general
only exists on a blowing-up $Y_k$ of
$Y$, and $Y_k$ will depend on $k$. While it  may be too much to
hope for a result for an arbitrary big bundle $B$, Demailly's
work  suggests  that  (7.3) should extend to all dimensions at
least for the sort of bundle occurring in the proof of Theorem
7.4.  See \cite{ELN} for a slightly different approach. \enddemo
\bl
 We will apply Proposition 7.3 in the spirit of \cite{De1},
\S8, to prove:
\proclaim{Theorem 7.4} Let $X$ be a smooth projective surface, $x
\in X$ a fixed point, and $s \ge 0$ a non-negative integer.
Suppose that $L$ is a nef line bundle on $X$ such that
$$L^2 \ge (s+2)^2 + 1 \ \text{ and } \ L \cdot C \ge s^2 + 3s +
3$$ for all curves $C  \subset X$ passing through $x$.  Then for
$k \gg 0$ there exists a divisor $D \in |kL|$ having an almost
isolated singularity of index $> (s+2)$ at
$x$. \endproclaim
\bl
 Propositon 7.2 then implies:
\proclaim{Corollary 7.5} If $L$ is a nef line bundle such that
$L^2 \ge 5$ and $L \cdot C \ge 3$ for all curves $C
\owns x$, then $\O_X(K_X + L)$ has a section which doesn't vanish
at $x$. More generally, if $L$ satisfies the inequalities in
Theorem 7.4, then $|K_X + L|$ generates $s$-jets at
$x$. \qed
\endproclaim
\bl
\ni See Exercise 7.6  for an analogous criterion for
$\O_X(K_X + L)$ to be very ample.  Note that the numbers are a
little weaker than the optimal bounds (Corollary 2.6) coming from
Reider's theorem.
\bl
\demo{Proof of Theorem 7.4} To simplify the exposition, we first
show how to produce a divisor $D
\in |kL|$ with an almost isolated singularity of index $\ge
(s+2)$ at $x$. At the end of the argument we will indicate the
(routine) changes necessary to get index $> (s+2)$.
\bl Let $f : Y \lra X$ be the blowing up of $X$ at $x$, with
$E \subset Y$ the exceptional divisor, and put $B = f^*L -
(s+2)E$.  Then
$$H^0(Y , \O_Y(kB)) = H^0(X, \O_X(kL) \otimes
\I_x^{k(s+2)}),$$ so by the usual dimension count (cf. proof of
Proposition 6.10):
$$
\aligned h^0(Y, \O_Y(kB)) & \ge \frac{k^2}{2} L^2 -
\binom{k(s+2) + 1}{2} + o(k^2) \\
	&=  \left ( L^2 - (s+2)^2 \right ) \frac{k^2}{2} + o(k^2).
\endaligned
$$ Consider as above the decomposition $|kB| = |M_k| + F_k$ of
$|kB|$ into its moving and fixed parts. Proposition 7.3 yields
the inequality
$$ M_k ^2 \ge \left ( L^2 - (s+2)^2 \right )k^2 + o(k^2).
$$ Hence by Hodge Index:
$$
\aligned M_k \cdot f^* L &\ge \sqrt{  \left ( L^2 \right )
\left ( M_k^2 \right )  } \\ & \ge k \left ( \sqrt{\left ( L^2
\right ) \left ( L^2 - (s+2)^2 \right ) } \right ) + o(k).
\endaligned
$$ Recalling that $F_k = f^*(kL) - k(s+2)E - M_k$, it follows that
$$ F_k \cdot f^*L \le k  \left ( L^2 - \sqrt{\left ( L^2
\right ) \left ( L^2 - (s+2)^2 \right ) }
\right ) + o(k). \tag *
$$ Now if $f_s(x)$ denotes the function
$$f_s(x) = x - \sqrt{x \left ( x - (s+2)^2 \right ) },$$ one
finds that $f_s(x) < s^2 + 3s + 3$ when $x \ge (s+2)^2 + 1$.
Hence it follows from (*) that for
$k \gg 0$,
$$F_k \cdot f^*L < (s^2 + 3s + 3)k  \tag **$$ provided that
$L^2 \ge (s+2)^2 + 1$.
\bl On the other hand, consider the image $$\overline F_k = f_*
F_k \subset X$$ of $F_k$ in $X$. Then
$\overline F_k$ is the fixed divisor of the linear series
$|\O_X(kL) \otimes \I_x^{k(s+2)} |$. It follows from (**) that
for $k \gg 0$:
$$\overline F_k \cdot L < (s^2 + 3s + 3)k.$$ But by hypothesis,
every effective curve $F \subset X$ passing through $x$ has
$L$-degree $ \ge s^2 + 3s + 3$. In particular, if $F$ is a
component of $\overline F_k$  through $x$, then
$ord_F(\overline F_k) < k$. Therefore Bertini's theorem implies
that a general divisor $D \in  |\O_X(kL) \otimes \I_x^{k(s+2)} |$
has multiplicity $< k$ in a punctured neighborhood of $x$, i.e.
$D$ has an almost isolated singularity of index $\ge (s+2)$ at
$x$.

\bl It remains to show that we can arrange to produce a divisor
$D \in |kL|$ with an almost isolated singularity of index
strictly greater than $s+2$. Had we proven a variant of
Proposition 7.3 for
$\Q$-divisors, one would simply go through the argument just
given with  the divisor $f^*L - (s+2)E$ replaced by $f^*L -
(s+2+\epsilon)E$ for small rational $\epsilon > 0$. To avoid
$\Q$-divisors in (7.3) one can equivalently choose a large
integer $m \gg 0$ and work with
$$B = f^*(mL) - \left ( m(s+2) + 1 \right ) E.$$ We leave details
to the reader. \qed \enddemo

\bl
\xca{Exercise 7.6} ({\bf Criterion for Very Ample Adjoint
Bundles}.)  Let $L$ be an ample line bundle on a smooth
projective surface $X$ such that $L^2 \ge 10$ and $L \cdot C
\ge 7$ for all curves $C \subset X$. Show that then
$\O_X(K_X + L)$ is very ample. [Given distinct points $x, y
\in X$, construct a divisor $D \in |kL|$ $(k \gg 0)$ with almost
isolated singularities of index $> 2$ at $x$ and
$y$.]  \qed \endxca

\bl
\xca{Exercise 7.7} ({\bf Effective  Matsusaka Theorem on
Surfaces}, \cite{FdB2}). Using the results of \cite{De1}, Siu
\cite{Siu1} has recently obtained an effective version of
Matsusaka's ``Big Theorem" (cf. \cite{LM},
\cite{KolM}) on varieties of all dimensions.
However the constants that appear in his statement are
exponential. Inspired by Siu's general result, Fern\'andez del
Busto shows in \cite{FdB2} that on a surface one can use
Proposition 7.3 to obtain the following considerably sharper
bound (having a different shape than Siu's):
\proclaim{Theorem} Let $L$ be an ample line bundle on the smooth
projective surface $X$, and set
$$ a = L^2 \ \ , \ \ b = (K_X + 4L) \cdot L.$$ If
$$m > \frac{ (b+1)^2}{2a} - 1 \ \ \bigg (\text{resp. } m >
\frac{ (b+1)^2}{2a} + 1  \bigg ),$$ then $\O_X(mL)$ is globally
generated {\rom (}resp. $\O_X(mL)$ is very ample{\rom )}.
 \endproclaim
\bl (i). Let $F$ and $G$ be ample line bundles on $X$. Show that
for $k \gg 0$:
$$h^0 \left ( X, \O_X(k(F-G)) \right ) \ge \left ( F^2 - 2F
\cdot G \right ) \frac{k^2}{2} + o(k^2).$$ [In fact, if $F^2 >
2F\cdot G$, then $h^0(\O_X(k(F-G))) \ge (F-G)^2
\frac{k^2}{2} + o(k^2)$.] In particular, taking
$ F = (m+3)L$ and $G = K_X + 4L$, it follows that
$$ h^0( \O_X(kB_m)) \ge \rho(m) \frac{k^2}{2} + o(k^2),$$ where:
$$
\gathered B_m = (m-1)L - K_X,  \\
\rho(m) = (m+3)^2 L^2 - 2(m+3)(K_X + 4L) \cdot L.
\endgathered$$
[Siu proves more generally that if $F$ and $G$ are ample line
bundles on a smooth projective variety $V$ of dimension $n$, and
if $\rho = F^n - n F^{n-1}\cdot G > 0$, then $h^0(V, \O_V(k(F-G))
\ge \rho \frac{k^n}{n!} + o(k^n)$.]
\sbl (ii). Fix $x \in X$, and suppose that $m$ is large enough so
that one has the inequalities:
$$
\gathered
\rho(m) > 4 \\ L \cdot B_m - \sqrt{(\rho(m) - 4)(L^2)} < 1.
\endgathered
\tag *
$$ Show that then  for $k \gg 0$ there exists a divisor $D
\in |kB_m|$ with an almost isolated singularity of index
$\ge 2$ at $x$. Deduce that $$H^1(X, \O_X(mL) \otimes \I_x) =
0,$$ and conclude that $|mL|$ free at $x$. [For the existence of
$D$, invoke Proposition 7.3  and argue as in the proof of Theorem
7.4. Then consider the $\Q$-divisor:
$$ M = mL - \frac{1}{k}D - K_X \equiv L.$$
 Apply vanishing to $f^*M$   on the blowing-up $f : Y \lra X$ of
$X$ at $x$.]
\sbl (iii). Show that the inequality $m > \frac{ (b+1)^2}{2a} -
1$ implies (*), which proves the first assertion of the Theorem.
Proceed similarly for very ampleness. \qed
\endxca

\bl
\xca{Exercise 7.8} ({\bf Bounding singularities via Seshadri
Constants.)} This exercise presents some further unpublished work
of {Geng Xu} showing how to deduce a variant of Theorem 7.4
from his bound (Exercise 5.16) on the Seshadri constant of an
ample line bundle along a finite set of general points.

\sbl
Let $L$ be an ample line bundle on a smooth projective surface
$X$ with $L^2 \ge 5$, and let $x \in X$ be a fixed point. The
aim is to show that if $L \cdot C \ge 5$ for all irreducible
curves $C \subset X$, then for $k \gg 0$ there exists a divisor $D
\in |kL|$ with an almost isolated singularity of index $> 2$ at
$x$.

\sbl
(i). Consider as above the decomposition
$$| \O_X(kL) \otimes  \I^{2k + 1}_x|  = \Sigma_k + F_k$$
of $| \O_X(kL) \otimes  \I^{2k + 1}_x|$ into its moving and
fixed parts. Given any $r = L^2 - 5$ points $x_1, \dots, x_r
\in X$, show that for $k \gg 0$ there exists a divisor $D^\prime
\in \Sigma_k$ with mult$_{x_i}(D^\prime) > k$ for all $1 \le i \le
r$.
\sbl
(ii). In the situation of (i), use Exercise 5.16 to prove
that if one chooses the $x_i$ sufficiently generally, then
$$L \cdot D^\prime > kr = k(L^2 - 5).$$
Deduce that $L \cdot F_k < 5k$ for $k \gg 0$, and conclude that
a general divisor $D \in | \O_X(kL) \otimes  \I^{2k + 1}_x|$ has
an almost isolated singularity of index $> 2$ at $x$. \qed
\endxca

\bl

\bl
\bl

\settabs\+University of California at Los Angeles and now is the
time \cr
\+ Robert LAZARSFELD \cr
\+ Department of Mathematics \cr
\+ University of California, Los Angeles \cr
\+ Los Angeles, CA  90024 \cr
\+ e-mail: \ \  rkl$\@$math.ucla.edu \cr

\end